\documentclass[12pt]{amsart}
\usepackage{amsmath, amsthm, amssymb}
\usepackage{bbm}
\usepackage[top=1in, bottom=1in, left=0.8in, right=0.8in]{geometry}
\usepackage{color}
\usepackage{hyperref}
\usepackage{indentfirst}
\usepackage{tikz}
\usepackage[tracking=true]{changes}

\hypersetup{hypertex=true,
	colorlinks=true,
	linkcolor=black,
    anchorcolor=black,
	citecolor=black}
\newtheorem{Thm}{Theorem}[section]
\newtheorem{Lm}[Thm]{Lemma}
\newtheorem{Prop}[Thm]{Proposition}
\newtheorem{Rmk}{Remark}
\numberwithin{equation}{section}
\newtheorem{Cor}[Thm]{Corollary}

\def\C{\mathbb C}

\def\R{\mathbb R}

\def\N{\mathbb N}

\renewcommand{\leq}{\leqslant}
\renewcommand{\geq}{\geqslant}

\def\p{\partial}

\DeclareMathOperator{\supp}{supp}

\begin{document}
	\title[Inverse problem for fractional equations on manifolds]{Inverse problem for fractional Schr\"{o}dinger equations with drift on closed Riemannian manifolds}

\author[T. Cai]{Tianyu Cai}
\address{Shanghai Center for Mathematical Sciences, Fudan University, Shanghai 200438, China. }
\email{23110840001@m.fudan.edu.cn}
\author[X. Chen]{Xi Chen}
\address{Shanghai Center for Mathematical Sciences, Fudan University, Shanghai 200438, China;	
	Center for Applied Mathematics, Fudan University, Shanghai 200433, China. }
\email{xi\_chen@fudan.edu.cn}
	\maketitle
	\begin{abstract} This paper is concerned about the inverse coefficient problems of variable-coefficient fractional Schr\"{o}dinger equations with drift on connected closed Riemannian manifolds.
	We prove that the knowledge of the underlying equation of order $\alpha\in (\frac{1}{2},1)$ on any non-empty open subset of the underlying manifold determines the Riemannian metric, the drift and the potential, simultaneously and uniquely, up to a gauge transformation, under the same geometric assumptions on the observation set as in \cite{feizmohammadi2024calderonproblemfractionalschrodinger}. The method of proof is based on that of \cite{feizmohammadi2024calderonproblemfractionalschrodinger} for fractional Schr\"{o}dinger operators, with the incorporation of the Runge approximation to recover the drift term.  
\end{abstract}

\section{Introduction}
A connected closed Riemannian manifold of dimension $n\geq 2$ is a connected compact boundaryless smooth $n$-manifold $M$ furnished with a Riemannian metric $g$. The gradient $\nabla_{g}$ and the Laplace--Beltrami operator $\Delta_{g}$, with respect to $g$, are locally given by
\begin{align*}
	\nabla_{g} u &= g^{jk} \frac{\p u}{\p x^j} \frac{\p }{\p x^k}, \\
	\Delta_{g} u &= |g|^{-\frac{1}{2}}\frac{\partial}{\partial x^{j}}\left(|g|^{\frac{1}{2}}g^{jk}\frac{\partial u}{\partial x^{k}}\right),
\end{align*}
where $g:=\left(g_{jk}\right)$, $|g|:=\det\left(g_{jk}\right)$, and $g^{-1}:=\left(g^{jk}\right).$

For $\alpha\in (\frac{1}{2},1)$, the fractional Schr\"{o}dinger operators with drift, denoted by $L_{ g,b,V}$, reads
\begin{equation}\label{MO}
	L_{ g,b,V}\left(u\right) := \left(-\Delta_{g}\right)^{\alpha}u + b\left(\nabla_{g} u\right) +Vu
\end{equation}
where the drift $b \in C^\infty\left(M, \mathbb{C} \otimes T^\ast M\right)$ is a complex-valued covector field,  and the potential $V \in C^\infty\left(M, \C\right)$ is a  complex-valued function. 

The formal $L^{2}$-adjoint operator $L_{ g,b,V}^{*}$ to $L_{ g,b,V}$ is defined through
\begin{equation*}
\left(u,L_{ g,b,V}^{*}\left(v\right)\right)_{L^{2}\left(M\right)} :=	\left(L_{ g,b,V}\left(u\right),v\right)_{L^{2}\left(M\right)}, \quad \forall u,v\in C^{\infty}\left(M\right).
\end{equation*}
More explicitly,  $L_{ g,b,V}^{*}$ is of the form
\begin{equation*}
	L_{ g,b,V}^{*}\left(u\right)=\left(-\Delta_{g}\right)^{\alpha}u +\nabla_{g}^{*}\left(\bar{b}u\right) + \overline{V}u
\end{equation*}
where $\nabla_g^{*}$ denotes the formal $L^{2}$-adjoint to $\nabla_g$. 

We are concerned about the fractional Schr\"{o}dinger equation with drift on $M$, 
\begin{equation}\label{FL_E_PP}
	L_{ g,b,V}\left(u\right) = f \,  
\end{equation}
where $f$ is supported in an open set $O\subset M$ and lies in the space
\begin{equation*}\label{SFS}
	H_{ g,b,V}^{O}:=\left\{f\in C_{c}^{\infty}\left(O\right):\left(f,v\right)_{L^{2}\left(M\right)}=0,\,\forall v\in \text{Ker}\left(L_{ g,b,V}^{*}\right) \right\}.
\end{equation*}
Moreover, it also involves the adjoint equation of \eqref{FL_E_PP}  \begin{equation}\label{FL_E_PP_adjoint}
	L_{ g,b,V}^\ast\left(u\right) = h \,  
\end{equation} where $h$ is an element of
	\begin{equation*}
	H_{ g,b,V}^{O,*} := \left\{h\in C_{c}^{\infty}\left(O\right): \left(h,v\right)_{L^{2}\left(O\right)} = 0, \: \forall v\in \text{Ker}\left(L_{ g,b,V}\right) \right\}.
\end{equation*}

 The standard Fredholm argument, e.g. Proposition \ref{WP} in Section \ref{StS1}, shows that the problems \eqref{FL_E_PP} and \eqref{FL_E_PP_adjoint} both have a unique solution, modulo $\mathrm{Ker}(L_{g,b,V})$ respectively $\mathrm{Ker}(L_{g,b,V}^{*})$. Hence, \eqref{FL_E_PP} admits a well-defined local
source-to-solution map, which reads
\begin{align}\label{STSM} 
	\Lambda_{ g,b,V}^{O}: H_{ g,b,V}^{O}  &\longrightarrow C^{\infty}\left(O\right) \\ \notag f&\longmapsto u^{f}|_{O}.
\end{align}

The inverse problem we aim to solve in this paper is:  \[\mbox{ Does  $\Lambda_{ g,b,V}^{O}$ simultaneously and uniquely determine $\left(g, b, V\right)$? } \]

If the metric $g$ is known, we give a positive answer for inversion of lower order coefficients.
\begin{Thm}\label{MR1}
	 Suppose $\left(M, g\right)$ is a connected closed Riemannian $n$-manifold, and $O\subset M$ is a nonempty open subset. Consider covector fields $b_1, b_2 \in C^\infty\left(M, \C \otimes T^\ast M\right)$ and potential functions $V_1, V_2 \in C^\infty\left(M, \C\right)$ which all vanish on $O$. For any $\alpha \in \left(\frac{1}{2},1\right)$, if the source-to-solution maps $	\Lambda_{ g,b_{1},V_{1}}^{O}$ and $\Lambda_{ g,b_{2},V_{2}}^{O}$ coincide, then it must hold that $b_1 = b_2$ and $V_1 = V_2$. 
\end{Thm}

If it needs to retrieve $\left(g, b, V\right)$ simultaneously, the answer to this inverse problem is negative in the usual sense.
This is a consequence of the geometric invariance $\Lambda_{ g,b,V}^{O}$ obeys.

In fact, a diffeomorphism $\Psi:M\rightarrow M$ induces the pull-back maps on Riemannian metrics, covector fields, and potential functions respectively. These maps are given by
\begin{equation*}
	\begin{aligned}
		\Psi^{*}g\left(x\right)&=\left(D\Psi\left(x\right)\right)^{T}g\left(\Psi\left(x\right)\right)D\Psi\left(x\right), && g \in C^\infty\left(M, T^\ast M \otimes T^\ast M\right); \\
		\Psi^{*}b\left(x\right)&=\left(D\Psi\left(x\right)\right)^{T}b\left(\Psi\left(x\right)\right), && b \in C^\infty\left(M, \C \otimes T^\ast M\right); \\
		\Psi^{*}V\left(x\right)&=V\left(\Psi\left(x\right)\right), && V \in C^\infty\left(M, \C\right).
	\end{aligned}
\end{equation*} 
We view $\Psi^\ast$ as the gauge transform for the coefficients in \eqref{FL_E_PP}
\begin{align}\label{PBA}
	\Psi^\ast : 
	\left(g,b,V\right)\longmapsto  \left(\Psi^{*}g,\Psi^{*}b,\Psi^{*}V\right).
\end{align} 
Theorem \ref{DG} verifies the gauge invariance of  the source-to-solution map \eqref{STSM}, i.e. 
\begin{equation}\label{eqn : gauge invariance StS}
	\Lambda^O_{ \Psi^{*}g,\Psi^{*}b,\Psi^{*}V} = \Lambda^O_{ g,b,V}.
\end{equation} Such a property suggests that the best one can achieve is the unique recovery of the gauge class of $\left(g, b, V\right)$.

From this perspective, we are able to reconstruct all coefficients simultaneously and uniquely.
\begin{Thm}\label{MR}
 Let $M$ be a connected closed Riemannian $n$-manifold and $O\subset M$ a nonempty connected open set. Assume for $j = 1, 2$  that
	\begin{itemize}
		\item  $g_j \in C^\infty(M, T^\ast M \otimes T^\ast M)$ is a Riemannian metric such that $g_1$ and $g_2$ agree on $O$; 
		\item $b_j \in C^\infty(M, \C \otimes T^\ast M)$ is a covector field, which vanishes on $O$;
		\item $V_j \in C^\infty(M, \C)$ is a potential function, which vanishes on $O$;
		\item $(M \setminus \bar{O}, g_j)$ is non-trapping;
		\item $O$ satisfies that there exists $p\in O$ such that 	\begin{equation}\label{eqn : FKU condition}
				\mathcal{A}_{ g_j}\left(p\right):=	\left\{q\in M:\mathrm{dist}_{g_j}\left(p,q\right)=\max\limits_{p^{\prime}\in M} \left\{\mathrm{dist}_{g_j}\left(p,p^{\prime}\right)\right\}\right\}\subset O.
			\end{equation} 
	\end{itemize}
 For any $\alpha \in \left(\frac{1}{2},1\right)$, if the source-to-solution maps $\Lambda_{ g_{1},b_{1},V_{1}}^{O}$ and $\Lambda_{ g_{2},b_{2},V_{2}}^{O}$ agree,
 then there must exist a diffeomorphism $\Psi: M\rightarrow M$ such that  
 \begin{equation*}
\Psi|_O = \mathrm{Id} \quad	\mbox{and}\quad	\left(g_{1},b_{1},V_{1}\right) =  \left(\Psi^{*}g_{2},\Psi^{*}b_{2},\Psi^{*}V_{2}\right).
	\end{equation*}
\end{Thm}

\begin{figure}
	\centering
\begin{tikzpicture}[x=0.75pt,y=0.75pt,yscale=-1,xscale=1]

	
	\draw   (207.06,93.12) .. controls (223,57.06) and (286.64,50.25) .. (349.21,77.92) .. controls (411.78,105.58) and (449.58,157.23) .. (433.64,193.29) .. controls (417.69,229.35) and (354.05,236.16) .. (291.48,208.5) .. controls (228.91,180.84) and (191.12,129.18) .. (207.06,93.12) -- cycle ;
	\draw   (226.5,88) .. controls (216.5,79) and (240,66.5) .. (263.5,64) .. controls (287,61.5) and (367.5,75) .. (411.5,141) .. controls (455.5,207) and (402.5,205) .. (373.5,134) .. controls (344.5,63) and (236.5,97) .. (226.5,88) -- cycle ;
	\draw    (268,130) .. controls (308,100) and (354.5,143) .. (355.5,173) ;
	\draw    (262.5,124) .. controls (274.5,140) and (324.5,199) .. (364.5,169) ;
	
	\draw (259,167) node [anchor=north west][inner sep=0.75pt]   [align=left] {M};
	\draw (285,67) node [anchor=north west][inner sep=0.75pt]   [align=left] {O};
\end{tikzpicture}
	\caption{Condition \eqref{eqn : FKU condition}}
\label{F2}
\end{figure}

\begin{Rmk}
	The condition \eqref{eqn : FKU condition} on the observation set $O$ was introduced in \cite{feizmohammadi2024calderonproblemfractionalschrodinger}. We refer the readers to \cite{feizmohammadi2024calderonproblemfractionalschrodinger, feizmohammadi2025inversespectralproblemssparse} for concrete examples of $O$ satisfying \eqref{eqn : FKU condition} with small volume (See Figure \ref{F2}). 

\end{Rmk}

The tool for the recovery of the drift term is the following Runge approximation for fractional equations.
\begin{Thm}\label{DR}
	 Let $\left(M,g\right)$ be a connected closed Riemannian $n$-manifold. Let $U$ be a nonempty open set with Lipschitz boundary such that $b\in C^{\infty}\left(M, \C\otimes T^{*}M\right)$ and $V\in C^{\infty}\left(M,\C\right)$ both vanish in $U$. For any $\alpha\in \left(\frac{1}{2},1\right)$,
	we consider the sets
		\begin{align}
			R&:=\left\{ \mbox{$u \mathbbm{1}_{M\setminus U}$ : $L_{g,b,V}(u)\in H_{g,b,V}^{U}$ and $u\perp \mathrm{Ker}(L_{g,b,V})$ in $L^{2}(M)$ } \right\}, \label{DR1}\\
				R^{*}&:=\left\{ \mbox{$u\mathbbm{1}_{M\setminus U}$ : $L_{g,b,V}^{*}(u)\in H_{g,b,V}^{U,*}$ and $u\perp \mathrm{Ker}(L_{g,b,V}^{*})$ in $L^{2}(M)$} \right\}. \notag
		\end{align}
	Then both $R\, \text{and}\, R^{*}\, \text{are dense in}\, \widetilde{H}^{\alpha}\left(M\setminus\bar{U}\right).$
\end{Thm}
\begin{Rmk}
	The vanishing assumption for $b$ and $V$ in $U$ guarantees that the Runge approximation in Theorem \ref{DR} is valid when $\dim(\mathrm{Ker}(L_{g,b,V}))\geq 1$.
\end{Rmk} 

Inverse problems of fractional elliptic equations are viewed as the nonlocal analogues of the classical Calder\'on problem.   Ghosh--Salo--Uhlmann \cite{MR4078233} proved, in a bounded Euclidean domain, that the potential function in fractional Schr\"{o}dinger equations can be uniquely determined by the exterior Dirichlet-to-Neumann type information. 
Bhattacharyya--Ghosh--Uhlmann \cite{MR4237942} generalized this framework to the Schr\"odinger equation with a nonlocal perturbation. Ceki\'c--Lin--R\"uland \cite{MR4092686} investigated inverse problems for the Schr\"odinger equation with a first order drift. They used the Runge approximation in the Euclidean domain to achieve the unique and simultaneous determination of the drift and the potential  from exterior DN type measurements.

With regard to principal coefficients, Covi, Railo, and Zimmermann \cite{ MR4062967, covi2022globalinversefractionalconductivity, MR4620353} studied inverse  problems of fractional  conductivity equations for isotropic media, and proved that the exterior DN type information uniquely recovers the conductivity. 
 
Anisotropic inverse problems are more involved. A natural strategy is to employ the famous Caffarelli--Silvestre extension theorem and convert nonlocal problems to local problems. See e.g. \cite{Ghosh02122017, ghosh2021calderonproblemnonlocaloperators, covi2023reductionfractionalcalderonproblem, MR4383014}.
However, this approach is utterly dependent on the results of corresponding local problems. In fact, the anisotropic problems of local elliptic equations, in general, are more challenging and even intractable.

Alternatively, a more enlightening thought for inverse problems is to exploit the nonlocality of the underlying equations for inversion, rather than to treat it as a hurdle. Feizmohammadi--Ghosh--Krupchyk--Uhlmann \cite{feizmohammadi2021fractionalanisotropiccalderonproblem} investigated inverse problems of fractional Poisson equations (with $b=0$, $V=0$)  on closed Riemannian manifolds. By digging into functional calculi for fractional operators, they developed an elegant argument, converting the coefficient identification of fractional Poisson  equations to that of wave equations  via the unique continuation of the  heat equations. Specifically, 
 \cite{feizmohammadi2021fractionalanisotropiccalderonproblem} proved that the knowledge of $\Lambda_{ g,0,0}^{O}$ determines the local data of the heat semigroup $e^{t\Delta_{g}}|_{O\times O}$ and further retrieves the wave kernel $\frac{\sin(t\sqrt{-\Delta_{g}})}{\sqrt{-\Delta_{g}}}|_{O\times O}$ via the Kannai's transmutation, which finally leads to the recovery of metrics from the wave kernel information, solved by Helin--Lassas--Oksanen--Saksala \cite[Theorem 2]{MR3826551}. In addition, R\"{u}land \cite{MR4874337} provided an alternative argument for the recovery of heat and wave kernels.

Recently, Feizmohammadi--Krupchyk--Uhlmann \cite{feizmohammadi2024calderonproblemfractionalschrodinger} considered the inverse problem for fractional Schr\"{o}dinger operators on closed Riemannian manifolds (with $b=0$) with receiver $O$ satisfying  \eqref{eqn : FKU condition}. Due to the presence of $V$, the information of $\Lambda_{g,0,V}^{O}$, via the framework in \cite{feizmohammadi2021fractionalanisotropiccalderonproblem}, only delivers 
  \begin{equation}\label{RM1I}
 	\left\{e^{t\Delta_{g}}(-\Delta_{g})^{\alpha}u^{f}|_{O}: t>0,f\in H_{g,0,V}^{O}\right\},
 \end{equation}
 which does not readily give the wave kernel data required in \cite[Theorem 2]{MR3826551}. 
By leveraging the spectral representation argument in \cite{MR4745418} and the entanglement principle \cite[Theorem 1.8]{feizmohammadi2024calderonproblemfractionalschrodinger}, \cite[Proposition 4.1]{feizmohammadi2024calderonproblemfractionalschrodinger} retrieved, from heat semigroup data \eqref{RM1I}, the non-normalized spectral data
  \begin{equation}\label{RM2}
  	\left\{\left(\lambda_{k}, \left\{\phi_{k,l}|_{O}\right\}_{l=1}^{d_{k}}\right)\right\}_{k=0}^{\infty}
  \end{equation}
 where each $\phi_{k,l}$ is an eigenfunction of $-\Delta_{g}$ corresponding to the eigenvalue $\lambda_{k}$ with multiplicity $d_{k}$ and $\left\{\left\{\phi_{k,l}\right\}_{l=1}^{d_{k}}\right\}_{k=0}^{\infty}$ forms a Schauder basis of $L^{2}(M)$. 
By using the finite speed of propagation and Tataru's unique continuation principle of wave equations, \cite[Theorem 1.11]{feizmohammadi2024calderonproblemfractionalschrodinger} established a highly non-trivial result that the non-normalized spectral data \eqref{RM2} reconstructs the normalized spectral data 
\begin{equation}\label{SPD}
  	\left\{\left(\lambda_{k},\{\phi_{k,l}|_{O}\}_{l=1}^{d_{k}}\right) \right\}_{k=0}^{\infty} \mbox{ such that $\left\{\left\{\phi_{k,l}\right\}_{l=1}^{d_{k}}\right\}_{k=0}^{\infty}$ are orthonormal}.
  \end{equation}
This step exploits the geometric assumption \eqref{eqn : FKU condition} and the non-trappingness of $(M\setminus\bar{O},g)$. In the end,  $g$ can be recovered from \eqref{SPD} by using \cite[Corollary 2]{MR3826551}.  See also \cite{MR4793146, feizmohammadi2025fractionalanisotropiccalderonproblem, feizmohammadi2025inversespectralproblemssparse} for other relevant results.
 
The novelty of our results is two-fold. On one hand, Theorem \ref{MR1} for lower order terms extends the result in \cite{MR4092686} to closed Riemannian manifolds, replacing the DN type information by the source-to-solution type measurements. For this purpose, we develop Theorem \ref{DR}, which serves as the Runge approximation on Riemannian manifolds. The proof relies on the techniques from \cite[Lemma 3.10]{feizmohammadi2024calderonproblemfractionalschrodinger}. On the other hand, Theorem \ref{MR} incorporates a drift term into the inverse result of \cite{feizmohammadi2024calderonproblemfractionalschrodinger} under the same assumptions on $O$ and for $\alpha\in (\frac{1}{2},1)$, using the method of \cite{feizmohammadi2024calderonproblemfractionalschrodinger}, together with the Runge approximation approach  in \cite{MR4092686}.
 
\subsection*{Outline of the paper} This paper is organized as follows. In Section \ref{S2}, we review the basic analytic tools needed for inversion. Section \ref{StS1} is devoted to some preliminary discussions. Specifically, we show that the local source-to-solution map $\Lambda_{ g,b,V}^{O}$ is well-defined (Proposition \ref{WP}) and that $\Lambda_{ g,b,V}^{O}$ is gauge invariant (Theorem \ref{DG}). Next, we prove the Runge approximation (Theorem \ref{DR}) in Section \ref{SRA}.  Finally, Section \ref{S4} recovers lower order terms on a priori known Riemannian manifolds  (Theorem \ref{MR1}); in the meanwhile, Section \ref{S3} reconstructs the Riemannian metric (Theorem \ref{MR}).
  
    \section{Preliminaries}\label{S2}
	\subsection{The fractional Laplacian}
	The Laplace--Beltrami operator $-\Delta_{g}$ is a self-adjoint operator on $L^{2}\left(M\right)$ with domain $D\left(- \Delta_{g}\right) = H^{2}\left(M\right)$. By the compactness of $M$, the spectrum $\text{Spec}\left(-\Delta_{g}\right)$ consists only of discrete eigenvalues, say $$\text{Spec}\left(-\Delta_{g}\right) = \left\{\lambda_k : k \in \mathbb{N},   \lambda_{0} = 0, \lambda_{k} <  \lambda_{k+1}\right\}.$$ 
For each $k \in \mathbb{N}_0$, we denote by $d_{k}$ the multiplicity of $\lambda_{k}$ and by $$V_{\lambda_{k}}:=\text{Ker}\left(-\Delta_{g} - \lambda_{k}\right)$$ the eigenspace associated with $\lambda_k$. In particular, the closedness of $M$ implies that $d_{0}=1$ and $V_{\lambda_{0}}=\text{Span}\{1\}$. 

For each $V_{\lambda_{k}}$, let $\{\phi_{k,1}, \cdots , \phi_{k,d_{k}}\}$ be an $L^{2}$-orthonormal basis of $V_{\lambda_{k}}$ and $\pi_{k}$ the orthogonal projection onto $V_{\lambda_{k}}$. Explicitly, $\pi_k$  is given by
	\begin{align} \pi_k : L^{2}\left(M\right) &\longrightarrow V_{\lambda_{k}},  \label{DP}\\
		 f &\longmapsto \sum\limits_{l=1}^{d_{k}}\left(f,\phi_{k,l}\right)_{L^{2}\left(M\right)}\phi_{k,l}, \notag
	\end{align}
	where $\left(\cdot,\cdot\right)_{L^{2}\left(M\right)}$ is the $L^{2}$-inner product on $M$.
	
 We employ the spectral theorem for $-\Delta_{g}$ on $L^{2}\left(M\right)$ to define the fractional Laplacian. 	For $\alpha\in \left(0,1\right)$,  the fractional Laplacian $\left(-\Delta_{g}\right)^{\alpha}$ is defined to be 
	\begin{equation*}
		\left(-\Delta_{g}\right)^{\alpha}u :=\sum\limits_{k=0}^{\infty}\lambda_{k}^{\alpha}\pi_{k}u
	\end{equation*}
for any $u$ lying in $$D\left(\left(-\Delta_{g}\right)^{\alpha}\right) :=H^{2\alpha}(M)= \left\{u\in L^{2}\left(M\right): \sum\limits_{k=0}^{\infty}\lambda_{k}^{2\alpha} \| \pi_{k}u\|_{L^{2}\left(M\right)}^{2}<\infty \right\}.$$
Then the inverse operator $\left(-\Delta_{g}\right)^{-\alpha}$ is given by
  \begin{equation}\label{IFL}
  	\left(-\Delta_{g}\right)^{-\alpha}v=\sum\limits_{k=1}^{\infty}\lambda_{k}^{-\alpha}\pi_{k}v,\quad \forall v \in D\left((-\Delta_{g})^{-\alpha}\right):=\left\{v\in L^{2}\left(M\right) : \left(v,1\right)_{L^{2}\left(M\right)}=0\right\}.
  \end{equation}   

\subsection{The fractional Sobolev space}
The fractional Sobolev space $H^{\alpha}\left(M\right)$ and the homogeneous Sobolev space $\dot{H}^{\alpha}\left(M\right)$ are subspaces of $L^{2}\left(M\right)$  defined respectively via the following  norms  
	\begin{align*}
		\|u \|_{H^{\alpha}\left(M\right)} &:=\left(\left\|\left(-\Delta_{g}\right)^{\frac{\alpha}{2}}u \right\|_{L^{2}\left(M\right)}^{2} + \left\|u \right\|_{L^{2}\left(M\right)}^{2} \right)^{\frac{1}{2}}= \left(\sum\limits_{k=0}^{\infty}\lambda_{k}^{\alpha}\|\pi_{k}u \|_{L^{2}\left(M\right)}^{2}+ \|u \|_{L^{2}\left(M\right)}^{2}\right)^{\frac{1}{2}},
\\
		\|u \|_{\dot{H}^{\alpha}\left(M\right)} &:=\left\|\left(-\Delta_{g}\right)^{\frac{\alpha}{2}}u \right\|_{L^{2}\left(M\right)}= \left(\sum\limits_{k=0}^{\infty}\lambda_{k}^{\alpha}\|\pi_{k}u \|_{L^{2}\left(M\right)}^{2}\right)^{\frac{1}{2}}.
	\end{align*}
	
	For a nonempty open set $U\subset M$, the space
	\begin{equation*}
		H^{\alpha}\left(U\right):=  \left\{u|_{U}:u\in H^{\alpha}\left(M\right)\right\}
	\end{equation*}
is	equipped with the norm 
	\begin{equation*}
		\|u \|_{H^{\alpha}\left(U\right)} := \inf \left\{\|w \|_{H^{\alpha}\left(M\right)}:w\in H^{\alpha}\left(M\right),\, w|_{U}=u\right\};
	\end{equation*}
 the space $H^{-\alpha}(U)$ is dual to
\begin{equation*}
	\widetilde{H}^{\alpha}(U) := \mbox{closure of $C_{c}^{\infty}(U)$ in $H^{\alpha}(M)$}.
\end{equation*}
When $U$ has Lipschitz boundary, 
\begin{equation*}
	\widetilde{H}^{\alpha}(U) = H_{\bar{U}}^{\alpha}: = \left\{u\in H^{\alpha}(M): \supp(u)\subset \bar{U}\right\}.
\end{equation*} 
\subsection{The source function set}	
  We next elucidate the structure of source functions for \eqref{FL_E_PP}.
	\begin{Lm}\label{RM}
		Under the same assumption in Theorem \ref{DR}, the source function set $H_{ g,b,V}^{U}$ satisfies the following properties. 
		\begin{itemize}
			\item When $\dim\left(\mathrm{Ker}\left(L_{ g,b,V}\right)\right)=0,$ then $$H_{ g,b,V}^{U} = C_{c}^{\infty}\left(U\right).$$
			
			\item When $\dim\left(\mathrm{Ker}\left(L_{ g,b,V}\right)\right)=N \geq 1$, there exists $\{\theta_{k}\}_{k=1}^{N}\subset C_{c}^{\infty}\left(U\right)$ such that $$f-\sum\limits_{k=1}^{N}\left(f,\eta_{k}\right)_{L^{2}\left(M\right)}\theta_{k}\in H_{ g,b,V}^{U},\quad \forall f\in C_{c}^{\infty}(U)$$ where $\{\eta_{k}\}_{k=1}^{N}$ is an $L^{2}\left(M\right)$-orthonormal basis of $\mathrm{Ker}\left(L_{ g,b,V}^{*}\right).$
		\end{itemize}
	\end{Lm}			
\begin{proof}
	See \cite[Remark 2.6 and Lemma 3.9]{feizmohammadi2024calderonproblemfractionalschrodinger}.
\end{proof}

\subsection{The heat kernel}
	The inversion will be achieved through the heat equation, which requires some classical results of the heat kernel. The first is the pointwise upper Gaussian type estimates on the heat kernel on manifolds from \cite{MR834612}.
	\begin{Lm}[Li--Yau estimates]\label{Li-Yau es}
		Let $(M,g)$ be a closed Riemannian $n$-manifold. Then there exist constants $C,c>0$ such that
	   	\begin{equation*}
	   	\left|e^{t\Delta_{g}}(x,y)\right| \leq Ct^{-\frac{n}{2}}e^{-c\frac{(\mathrm{dist}_{g}(x,y))^{2}}{t}}, \quad t\in (0,1), \: x,y\in M.
	   \end{equation*}
	\end{Lm}
	
	We also need Gaussian type integral estimates for heat kernel due to \cite{MR102097, MR1226938}.
	\begin{Lm}[Davies--Gaffney inequalities]\label{LM_DG}
	Let $M$ be a complete Riemannian manifold. If A and B are two disjoint measurable subsets and $f,h\in L^{2}\left(M\right)$ with $\mathrm{supp}\left(f\right)\subset A$, $\mathrm{supp}\left(h\right)\subset B$, then 
		\begin{equation}\label{DGI}
			\left|\left(e^{t\Delta_{g}}f,h\right)_{L^{2}\left(M\right)}\right| \leq \|f \|_{L^{2}\left(A\right)}\|h \|_{L^{2}\left(B\right)}e^{-\frac{\mathrm{dist}_{g}^{2}\left(A,B\right)}{4t}-\lambda_{0}t},
			\end{equation}
			where $\lambda_{0}$ is the minimum eigenvalue of $-\Delta_{g}.$
		\end{Lm}
	The Gaussian decay in Davies--Gaffney inequalities \eqref{DGI} is used through the following lemma to prove the uniqueness of some heat kernel data.
	\begin{Lm}\label{DG_E}
		Let $\left(M,g\right)$ be a closed Riemannian $n$-manifold and $O\subset M$ a nonempty open set. Suppose $b\in C^{\infty}\left(M,\C \otimes T^{*}M\right)$ and $V\in C^{\infty}\left(M, \C\right)$ both vanish on $O$. There holds that
		\begin{equation}
			\left\|e^{\frac{1}{s}\Delta_{g}}\left(f-b\left(\nabla_{g} u\right)-Vu\right) \right\|_{L^{1}\left(\omega\right)} \leq Ce^{-cs}
			\end{equation}
		where $K := \mathrm{supp}\left(f\right) \subset O$, $\omega\Subset O$ with $\overline{\omega}\cap K=\emptyset$. Here $C>0$ is a constant depending on $M,\, g, \, b,\, V,\, O,\, f $ and $c>0$ is a constant depending only on $\mathrm{dist}_{g}\left(K\cup  \left(M\setminus\bar{O}\right), \omega\right)$.
	\end{Lm}
	\begin{proof}
		Since $K\cup \left(M\setminus\bar{O}\right)$ and $\omega$ are disjoint, $$\text{dist}_{g}\left(K\cup \left(M\setminus\bar{O}\right), \omega\right)>0.$$ By \eqref{DGI} and $\lambda_{0}=0$, there holds for any $h\in C_{c}^{\infty}\left(\omega\right)$ that 
		\begin{equation*}
			\begin{aligned}
				&\left|\left(e^{\frac{1}{s}\Delta_{g}}\left(f-b\left(\nabla_{g} u\right)-Vu\right),h\right)_{L^{2}\left(M\right)}\right|\\
				&\leq \left\|f-b\left(\nabla_{g} u\right)-Vu \right\|_{L^{2}\left(K\cup \left(M\setminus\bar{O}\right) \right)} \left\|h\right\|_{L^{2}\left(\omega\right)} e^{-s\frac{\text{dist}_{g}^{2}\left(K\cup \left(M\setminus\bar{O}\right),\omega\right)}{4}}. \\
			\end{aligned}
		\end{equation*}
		It implies that 
		\begin{equation*}\label{LM4a}
			\left\|e^{\frac{1}{s}\Delta_{g}}\left(f-b\left(\nabla_{g} u\right)-Vu\right) \right\|_{L^{2}\left(\omega\right)} \leq \left\|f-b\left(\nabla_{g} u\right)-Vu \right\|_{L^{2}\left(K\cup \left(M\setminus\bar{O}\right)\right)}e^{-s\frac{\text{dist}_{g}^{2}\left(K\cup \left(M\setminus\bar{O}\right),\omega\right)}{4}}.
		\end{equation*}
	Since $\omega$ is compact, we obtain that
		\begin{equation*}
			\left\|e^{\frac{1}{s}\Delta_{g}}\left(f-b\left(\nabla_{g} u\right)-Vu\right) \right\|_{L^{1}\left(\omega\right)} \leq Ce^{-cs}.
		\end{equation*}
	\end{proof}

The unique continuation principle of heat equations, due to \cite{MR1024191}, is vital for the recovery of heat semigroup data.
	\begin{Lm}[Unique continuation of heat equations]\label{LUCP}
	Let $\Omega$ be an open, bounded, connected subset of $\mathbb{R}^{n}$. Let $u\in L^{2}\left(\left(0,T\right), H^{1}(\Omega)\right)$ be a solution of 
	\begin{equation*}
		\partial_{t}u=\Delta u \quad \text{in} \quad \Omega\times \left(0,T\right).
	\end{equation*}
	If in addition there exists some $\left(x_{0},t_{0}\right)\in \Omega\times \left(0,T\right)$ such that $$\int_{B_{r}\left(x_{0}\right)}u^{2}\left(x,t_{0}\right)dx=o\left(r^{N}\right)\quad \mbox{as $r\to 0^{+}$} $$  for all integer $N$, then $u\left(x,t_{0}\right)\equiv 0$ on $\Omega$.
\end{Lm}

\subsection{Inverse spectral problems}
The metric will be recovered from the spectral data. We need some results concerning inverse spectral problems.
	    \begin{Lm}\label{SK}
		Assume $\alpha\in \left(\frac{1}{2},1\right).$ Let $\left(M,g\right)$ be a connected closed Riemannian manifold and $O\subset M$ a nonempty open set. Suppose $b\in C^{\infty}\left(M,\mathbb{C} \otimes T^{*}M\right)$ and $V\in C^{\infty}\left(M, \mathbb{C}\right)$ both vanish on $O$. If $\phi\in C^{\infty}\left(M\right)$ is an eigenfunction of $-\Delta_{g}$ corresponding to some eigenvalue $\lambda>0$, then there exists $$u\in \{u : \mbox{$L_{g,b,V}(u)\in H_{ g,b,V}^{O}$ and $u\perp \mathrm{Ker}(L_{g,b,V})$ in $L^{2}(M)$}\}$$ such that $\left(u,\phi\right)_{L^{2}\left(M\right)}\neq 0.$
	\end{Lm}
	\begin{Rmk}
		The proof of Lemma \ref{SK} is similar to  \cite[Lemma 3.10]{feizmohammadi2024calderonproblemfractionalschrodinger}.
		\end{Rmk}

It was proved in  \cite{MR3826551} that the local spectral data recovers the metric up to a diffeomorphism. 
	\begin{Lm}\label{ML}
		Let $\left(M,g\right)$ be a connected closed Riemannian $n$-manifold, and $\mathcal{X}$ a nonempty open subset. Assume that  $\{\varphi_{k}\}_{k=1}^{\infty}\subset C^{\infty}\left(M\right)$ is the collection of orthonormal eigenfunctions of $\Delta_{g}$ in $L^{2}\left(M\right)$ associated with eigenvalues $\{\lambda_{k}\}_{k=1}^{\infty}$. Then the spectral data 
		\begin{equation*}
			\{\mathcal{X}, \{\varphi_{k}|_{\mathcal{X}}\}_{k=1}^{\infty}, \{\lambda_{k}\}_{k=1}^{\infty}\}
		\end{equation*}
		determines $\left(M,g\right)$ up to a diffeomorphism.
	\end{Lm}

	\section{The source-to-solution map}\label{StS1}
		We begin with the weak solution for equation \eqref{FL_E_PP}. For $f\in H^{-\alpha}\left(M\right)$, we say that $u\in H^{\alpha}\left(M\right)$ is a weak solution to \eqref{FL_E_PP} if
		\begin{equation}\label{SES}
			\left(L_{g,b,V}\left(u\right),v\right)_{L^{2}\left(M\right)} = \left(f,v\right)_{H^{-\alpha}\left(M\right),H^{\alpha}\left(M\right)},\, \forall v\in C^{\infty}\left(M\right).
		\end{equation}
	We justify the well-definedness of $\Lambda_{ g,b,V}^{O}$ defined in \eqref{STSM}. This is a consequence of the following solvability results of equations \eqref{FL_E_PP} and its adjoint equation.  	  
	\begin{Prop}\label{WP}
		Let $ M ,g, b, V, \alpha$ be as in Theorem \ref{DR}. Then there holds
 $$\dim\left(\mathrm{Ker}\left(L_{ g,b,V}\right)\right)=\dim\left(\mathrm{Ker}\left(L_{ g,b,V}^{*}\right)\right)<\infty;$$ 
		 for given $f\in H^{-\alpha}\left(M\right)$, the inhomogeneous  equation
			\eqref{FL_E_PP}
			is solvable if and only if $$\left(f,v\right)_{H^{-\alpha}(M),H^{\alpha}(M)}=0,\, \forall v\in \mathrm{Ker}\left(L_{ g,b,V}^{*}\right);$$	 
		for given $h\in H^{-\alpha}\left(M\right)$,	the inhomogeneous  equation  
				\begin{equation*}\label{WP_AE}
			L_{ g,b,V}^{*}\left(u\right)= h
		\end{equation*}
			is solvable if and only if $$\left(h,v\right)_{H^{-\alpha}(M),H^{\alpha}(M)}=0,\, \forall v\in \mathrm{Ker}\left(L_{ g,b,V}\right).$$	 
	\end{Prop}
	\begin{proof}		
	By the standard Fredholm alternative argument, it suffices to show that $L_{ g,b,V}$ is a Fredholm operator. The proof works verbatim for $L_{ g,b,V}^\ast$. 
	
	For this purpose, we prove instead that there exists a constant $C_{  g, b} > 0$ such that $$L_{ g,b}:=\left(-\Delta_{g}\right)^{\alpha}+ b\left(\nabla_{g} \cdot\right) +C_{  g, b}$$ is invertible, where $C_{  g, b}$ is given explicitly in \eqref{eqn : C(Mgb)}  at the end of the proof. Then the desired Fredholm property of $L_{ g,b,V}$ follows from the compactness of multiplication by $V-C_{  g, b}$.	
				
 To prove the invertibility of $L_{  g, b}$,	we introduce the following sesquilinear form $B_{ g,b}\left(u,v\right)$ on $H^\alpha\left(M\right)$ associated with $L_{ g,b}$, which is explicitly given by
	\begin{equation*}B_{ g,b}\left(u,v\right) = \left(L_{ g,b}(u),v\right)_{L^{2}(M)}.
	\end{equation*}			
	It is clear that $B_{ g,b}\left(u,v\right)$ is linear in $u$ and skewlinear in $v$.
		We shall show that $B_{ g,b}$ is bounded and coercive on $H^\alpha\left(M\right)$.
		Then the Lax--Milgram Lemma guarantees that for any  $\mathfrak{w} \in H^{-\alpha}\left(M\right)$ there exists a unique $\mathfrak{u}\in H^{\alpha}(M)$ such that 
		\[B_{ g,b}\left(\mathfrak{u},      v\right) =  \left(\mathfrak{w}, v\right)_{H^{-\alpha}(M),H^{\alpha}(M)}, \quad \forall v \in H^\alpha\left(M\right).\] Consequently, the inverse $L^{-1}_{  g, b}$ is well-defined 
		\begin{align*}L^{-1}_{  g, b} : H^{-\alpha}\left(M\right) &\longrightarrow H^{\alpha}\left(M\right)  \\ \mathfrak{w} &\longmapsto \mathfrak{u} \end{align*}
		
		  The boundedness of $B_{ g,b}\left(\cdot,\cdot\right)$ on $H^{\alpha}\left(M\right)\times H^{\alpha}\left(M\right)$ is proved as follows.
			\begin{equation*}\label{WP_1}
				\begin{aligned}
	    	\left|B_{ g,b}\left(u,v\right)\right| &\leq \left|\left(\left(-\Delta_{g}\right)^{\frac{\alpha}{2}}u,\left(-\Delta_{g}\right)^{\frac{\alpha}{2}}v\right)_{L^{2}\left(M\right)}\right| + 
	    	\left|\left(b\left(\nabla_{g} u\right),v\right)_{L^{2}\left(M\right)}\right| + C_{  g, b}\left|\left(u,v\right)_{L^{2}\left(M\right)}\right| \\
	    	& \leq \| u\|_{\dot{H}^{\alpha}\left(M\right)} \| v\|_{\dot{H}^{\alpha}\left(M\right)} + \left|\left(b\left(\nabla_{g} u\right),v\right)_{L^{2}\left(M\right)}\right| + C_{  g, b}\|u \|_{L^{2}\left(M\right)} \|v\|_{L^{2}\left(M\right)}\\
		   & \leq \left(1+C_{  g, b}\right)\|u\|_{H^{\alpha}\left(M\right)} \|v \|_{H^{\alpha}\left(M\right)} + \left|\left(b\left(\nabla_{g} u\right),v\right)_{L^{2}\left(M\right)}\right|.
		\end{aligned}
				\end{equation*}
		Here the second term on the RHS obeys the following estimates
				\begin{equation*}\label{WP_2}
						\begin{aligned}
						\left|\left(b\left(\nabla_{g} u\right),v\right)_{L^{2}\left(M\right)}\right| &  \leq \left\|v\bar{b} \right\|_{H^{1-\alpha}\left(M\right)} \left\| \nabla_{g} u \right\|_{H^{\alpha-1}\left(M\right)} \\
						&\leq \left(\|b\|_{L^{\infty}\left(M\right)}\|v \|_{H^{1-\alpha}\left(M\right)}+\left\|\left(1-\Delta_{g}\right)^{\frac{1-\alpha}{2}}b \right\|_{L^{\infty}\left(M\right)}\|v \|_{L^{2}\left(M\right)} \right) \left\| \nabla_{g} u \right\|_{H^{\alpha-1}\left(M\right)} \\
						&\leq \left(\|b\|_{L^{\infty}\left(M\right)}+ \left\|\left(1-\Delta_{g}\right)^{\frac{1-\alpha}{2}}b \right\|_{L^{\infty}\left(M\right)}\right) \| v\|_{H^{1-\alpha}\left(M\right)}\sqrt{n} \|u \|_{H^{\alpha}\left(M\right)}.  
					\end{aligned}
					\end{equation*}
			Noting $0 < 1-\alpha < \frac{1}{2} <\alpha < 1$, we have	\begin{equation*}\label{WP_2}
					\begin{aligned}
						\left|\left(b\left(\nabla_{g} u\right),v\right)_{L^{2}\left(M\right)}\right| 
						&\leq \left(\|b\|_{L^{\infty}\left(M\right)}+\left\|\left(1-\Delta_{g}\right)^{\frac{1-\alpha}{2}}b \right\|_{L^{\infty}\left(M\right)}\right)\sqrt{n} \| v\|_{H^{\alpha}\left(M\right)} \|u \|_{H^{\alpha}\left(M\right)}.
					\end{aligned}
				\end{equation*} 			
				
		With regard to	the	coercivity of $B_{ g,b}$, we observe that $B_{ g,b}\left(v,v\right)$ is bounded from below by
			\begin{equation}\label{WP_LM}
				\begin{aligned}
					B_{ g,b}\left(v,v\right) 
					&\geq \left\|\left(-\Delta_{g}\right)^{\frac{\alpha}{2}}v \right\|_{L^{2}\left(M\right)}^{2} -  \left|\int_{M}\bar{v}b\left(\nabla_{g} v\right)dx\right|+ C_{g,b}\|v \|_{L^{2}\left(M\right)}^{2}.
				\end{aligned}
			\end{equation}
			The second term on the RHS can be estimated via the interpolation inequality in  \cite[Theorem 5.2]{MR450957}. That is, for $\alpha\in \left(\frac{1}{2},1\right)$,
			\begin{equation}\label{II}
				\| v\|_{H^{1-\alpha}\left(M\right)} \leq C_1 \|v \|_{L^{2}\left(M\right)}^{\frac{2\alpha-1}{\alpha}} \|v \|_{H^{\alpha}\left(M\right)}^{\frac{1-\alpha}{\alpha}}.
			\end{equation} 
			It follows that
			\begin{equation*}
				\begin{aligned}
					\left|\int_{M}\bar{v}b\left(\nabla_{g} v\right)dx\right| &\leq \left\|v\bar{b} \right\|_{H^{1-\alpha}\left(M\right)} \left\| \nabla_{g} v \right\|_{H^{\alpha-1}\left(M\right)}  \\
					&\leq \left(\|b\|_{L^{\infty}\left(M\right)}+   \left\|\left(1-\Delta_{g}\right)^{\frac{1-\alpha}{2}}b \right\|_{L^{\infty}\left(M\right)}\right)\|v \|_{H^{1-\alpha}\left(M\right)}\sqrt{n} \|v \|_{H^{\alpha}\left(M\right)} \\
					& \leq C_1 \sqrt{n}\left(\|b\|_{L^{\infty}\left(M\right)}+ \left\|\left(1-\Delta_{g}\right)^{\frac{1-\alpha}{2}}b \right\|_{L^{\infty}\left(M\right)} \right)\|v \|_{L^{2}\left(M\right)}^{\frac{2\alpha-1}{\alpha}} \|v \|_{H^{\alpha}\left(M\right)}^{\frac{1}{\alpha}}  
				\end{aligned}
			\end{equation*} 	provided some independent  constant $C_{1}>0$ from \eqref{II}.
	 Moreover, applying	Young's inequality with $p=2\alpha>1$ and $q=\frac{2\alpha}{2\alpha-1}>1$,  
		\begin{equation*}
			ab \leq \frac{a^{p}}{p} + \frac{b^{q}}{q}, \quad \forall a,b\geq 0,
		\end{equation*}    yields that
				\begin{align}\label{WP_BE_3} 
				 \left|\int_{M}\bar{v}b\left(\nabla_{g} v\right)dx\right| &\leq  C_{2} \left(\varepsilon^{\frac{1}{2\alpha}}\|v \|_{H^{\alpha}\left(M\right)}^{\frac{1}{\alpha}}\right) \left(\varepsilon^{-\frac{1}{2\alpha}}\|v \|_{L^{2}\left(M\right)}^{\frac{2\alpha-1}{\alpha}}\right)\\
			\notag	& \leq C_{2}\max\left(\frac{1}{2\alpha},\frac{2\alpha-1}{2\alpha}\right)
				\left(\varepsilon\|v \|_{H^{\alpha}\left(M\right)}^{2}+ \varepsilon^{-\frac{1}{2\alpha-1}}\|v\|_{L^{2}\left(M\right)}^{2}\right). 
		\end{align}
		with $$C_{2}:= C_1\sqrt{n} \left(\|b\|_{L^{\infty}\left(M\right)}+\left\|\left(1-\Delta_{g}\right)^{\frac{1-\alpha}{2}}b \right\|_{L^{\infty}\left(M\right)}  \right).$$
Inserting \eqref{WP_BE_3} into \eqref{WP_LM} delivers that
			\begin{equation*}
				\begin{aligned}
				 {B_{ g,b}\left(v,v\right)} 
				 &\geq \left\|\left(-\Delta_{g}\right)^{\frac{\alpha}{2}}v \right\|_{L^{2}\left(M\right)}^{2} + \|v \|_{L^{2}\left(M\right)}^{2} + \left(C_{  g, b}-1\right)\|v \|_{L^{2}\left(M\right)}^{2} \\
					&\quad - C_{2}\max(\frac{1}{2\alpha},\frac{2\alpha-1}{2\alpha})
					\left(\varepsilon \|v \|_{H^{\alpha}\left(M\right)}^{2} + \varepsilon^{-\frac{1}{2\alpha-1}}\|v \|_{L^{2}\left(M\right)}^{2}  \right)  \\
					& = \|v\|_{H^{\alpha}\left(M\right)}^{2} - \varepsilon C_{2}\max(\frac{1}{2\alpha},\frac{2\alpha-1}{2\alpha})\|v \|_{H^{\alpha}\left(M\right)}^{2} + \\
					&\quad  \left(C_{  g, b}-1-\varepsilon^{-\frac{1}{2\alpha-1}}C_{2}\max(\frac{1}{2\alpha},\frac{2\alpha-1}{2\alpha})\right)\|v \|_{L^{2}\left(M\right)}^{2} 
				\end{aligned}
			\end{equation*}			 
	Finally,	the proof of	 coercivity
			 \[
			B_{ g,b}\left(v,v\right)	  \geq \frac{1}{2}\|v \|_{H^{\alpha}\left(M\right)}^{2}\] is concluded by letting \begin{align}  \label{eqn : C(Mgb)} \left\{\begin{aligned}\varepsilon&:=\frac{1}{2C_{1}\max(\frac{1}{2\alpha},\frac{2\alpha-1}{2\alpha})\sqrt{n}\left(\|b\|_{L^{\infty}\left(M\right)}+\|\left(1-\Delta_{g}\right)^{\frac{1-\alpha}{2}}b \|_{L^{\infty}\left(M\right)}\right)}\\ 
				C_{  g, b}&:=2+ \varepsilon^{-\frac{1}{2\alpha-1}}C_{1}\max(\frac{1}{2\alpha},\frac{2\alpha-1}{2\alpha})\sqrt{n}\left(\|b\|_{L^{\infty}\left(M\right)}+ \left\|\left(1-\Delta_{g}\right)^{\frac{1-\alpha}{2}}b \right\|_{L^{\infty}\left(M\right)}\right). \end{aligned} \right.\end{align} where $C_{1}>0$ is the constant in \eqref{II}.
		\end{proof}		
		Using Proposition \ref{WP}, we can define the local source-to-solution map $\Lambda_{ g,b,V}^{O}$. 
		\begin{Cor}\label{well-defined StS}
		For $f\in H_{ g,b,V}^{O}$,  there exists a unique solution $u^{f}\in C^{\infty}\left(M\right)$ to \eqref{FL_E_PP} subject to the condition $\left(u^{f},v\right)_{L^{2}\left(M\right)} =0$ for all $v\in \mathrm{Ker}\left(L_{ g,b,V}\right)$. Consequently, the local source-to-solution map
		\begin{equation*}\label{LSTS}
			 \Lambda_{ g,b,V}^{O}\left(f\right) :=u^{f}|_{O},\quad \forall f\in H_{ g,b,V}^{O},
	   \end{equation*} is well-defined.
   \end{Cor}
	
    We conclude this section by proving the gauge invariance of the local source-to-solution map.
    		\begin{Thm}\label{DG}
	Let $ M, O, g, b, V, \alpha$ be as in Theorem \ref{DR}.
Suppose $\Psi:M\rightarrow M$ is a diffeomorphism such that $\Psi|_{O}=\mathrm{Id}$. Then $\Lambda_{g, b, V}^{O}$ is invariant under the action of $\Psi$-action on $\left(g, b, V\right)$ defined in \eqref{PBA}. Namely, \eqref{eqn : gauge invariance StS} holds.
\end{Thm}
\begin{proof}  
	We first prove that for $u\in C^{\infty}\left(M\right)$, there holds
	\begin{equation}\label{NGT}
		L_{\Psi^{*}g,\Psi^{*}b,\Psi^{*}V}\left(u\circ\Psi\right) = L_{g,b,V}\left(u\right)\circ\Psi 
	\end{equation}	
A standard computation of change of coordinate shows that  the principal terms satisfy  \begin{equation}\label{API1}\left(-\Delta_{\Psi^{*}g}\right)^{\alpha}\left(u\circ \Psi\right) =\left(\left(-\Delta_{g}\right)^{\alpha}u\right)\circ \Psi.\end{equation}
The following identity of  the pull-back of the zeroth order term is also clear. 
\begin{equation}\label{API2}
	\Psi^{*}\left(Vu\right)    =  \left(V\circ\Psi\right)\left(u\circ\Psi\right) = \left(Vu\right)\circ\Psi.
\end{equation}
Regarding the first order terms, one can prove 
\begin{equation}\label{API3}
	\Psi^{*}b\left(\nabla_{\Psi^{*}g} \left(u\circ\Psi\right)\right) = \left(b\left(\nabla_{g} u\right)\right)\circ\Psi.
\end{equation}
 by plugging in the gradient $\nabla_{\Psi^\ast g}$ and the covector field $\Psi^\ast b$.  
 More precisely, we invoke \begin{align*}
 	\Psi^{*}g\left(x\right)&=D\Psi\left(x\right)^{T}g\left(\Psi\left(x\right)\right)D\Psi\left(x\right), \\ \Psi^{*}b\left(x\right)&=D\Psi\left(x\right)^{T}b\left(\Psi\left(x\right)\right), 
 	\end{align*}
 	to calculate
\begin{align*}
	\lefteqn{ \Psi^{*}b\left(\nabla_{\Psi^{*}g}\left(u\circ\Psi\right)\right)\left(x\right) }\\
	&=\left(\Psi^{*}b\left(x\right)\right)^{T}\left(\Psi^{*}g\left(x\right)\right)^{-1}\left(D\Psi\left(x\right)\right)^{T}Du\left(\Psi\left(x\right)\right) \\
	& =	\left(\left(D\Psi\left(x\right)\right)^{T}b\left(\Psi\left(x\right)\right)\right)^{T}\left(\left(D\Psi\left(x\right)\right)^{T}g\left(\Psi\left(x\right)\right)D\Psi\left(x\right)\right)^{-1}\left(\left(D\Psi\left(x\right)\right)^{T}Du\left(\Psi\left(x\right)\right)\right)\\
	&= \left(b\left(\Psi\left(x\right)\right)\right)^{T}D\Psi\left(x\right)\left(D\Psi\left(x\right)\right)^{-1}\left(g\left(\Psi\left(x\right)\right)\right)^{-1}\left(D\Psi\left(x\right)^{T}\right)^{-1}\left(D\Psi\left(x\right)\right)^{T}Du\left(\Psi\left(x\right)\right) \\
	&=\left(b\left(\Psi\left(x\right)\right)\right)^{T}\left(g\left(\Psi\left(x\right)\right)\right)^{-1}Du\left(\Psi\left(x\right)\right) \\
	&= b\left(\nabla_{g} u\right)\left(\Psi\left(x\right)\right).
	\end{align*}
Therefore, \eqref{NGT} is a consequence of \eqref{API1}, \eqref{API2} and \eqref{API3}.

Next, we prove that 
\begin{equation*}
H_{\Psi^{*}g,\Psi^{*}b,\Psi^{*}V}^{O} = H_{g,b,V}^{O}
\end{equation*}
which, by definition, amounts to $$\text{Ker}\left(L_{ \Psi^{*}g,\Psi^{*}b,\Psi^{*}V}^{*}\right)|_{O} = \text{Ker}\left(L_{ g,b,V}^{*}\right)|_{O}.$$
In fact, we use \eqref{NGT} to compute for $u\in \text{Ker}\left(L_{ g,b,V}^{*}\right)$ that  
\begin{equation*}\label{GTT2}
\begin{aligned}
\left(L_{ \Psi^{*}g,\Psi^{*}b,\Psi^{*}V}^{*}\left(u\circ\Psi\right),v\circ\Psi\right)_{L^{2}_{ \Psi^{*}g}} &= \left(u\circ \Psi,L_{ \Psi^{*}g,\Psi^{*}b,\Psi^{*}V}\left(v\circ\Psi\right)\right)_{L^{2}_{ \Psi^{*}g}} \\
&=\left(u\circ\Psi,L_{ g,b,V}\left(v\right)\circ\Psi\right)_{L^{2}_{ \Psi^{*}g}} \\
&=\left(u,L_{ g,b,V}\left(v\right)\right)_{L^{2}_{g}} \\
&=\left(L_{ g,b,V}^{*}\left(u\right),v\right)_{L^{2}_{g}}=0, \quad \forall v\in C^{\infty}\left(M\right).
\end{aligned}
\end{equation*}
It follows that $$u\circ \Psi\in \text{Ker}\left(L_{ \Psi^{*}g,\Psi^{*}b,\Psi^{*}V}^{*}\right).$$ 
Noting $\Psi|_{O}=\mathrm{Id}$, we see that $$\text{Ker}\left(L_{ g,b,V}^{*}\right)|_{O}\subset \text{Ker}\left(L_{ \Psi^{*}g,\Psi^{*}b,\Psi^{*}V}^{*}\right)|_{O}.$$ The opposite inclusion can be proved verbatim.

Now we are ready to prove the gauge invariance \eqref{eqn : gauge invariance StS}.
Take $f\in H_{ g,b,V}^{O}$ and let $u^{f}$ solve
\begin{equation*}\label{GTT1}
 L_{ g,b,V}\left(u^{f}\right)=f  
\end{equation*} 
subject to $$(u^{f},v)_{L^{2}(M)}=0,\, \forall v\in \text{Ker}(L_{g,b,V}).$$
By \eqref{NGT},  $u^{f}\circ\Psi$ obeys
$$
L_{ \Psi^{*}g,\Psi^{*}b,\Psi^{*}V}\left(u^{f}\circ\Psi\right) = f\circ\Psi. 
$$
This means
\[\Lambda^O_{ \Psi^{*}g,\Psi^{*}b,\Psi^{*}V}\left(f\circ\Psi\right) = \left(u^{f}\circ\Psi\right)|_{O}=\Lambda_{ g,b,V}^{O}(f)\circ\Psi. \]
Combining with the fact $$\left(u^{f}\circ \Psi,v\circ\Psi\right)_{L^{2}_{\Psi^{*}g}\left(M\right)}= (u^{f},v)_{L^{2}_{g}(M)}=0,  \quad \forall v\circ\Psi\in \text{Ker}\left(L_{ \Psi^{*}g,\Psi^{*}b,\Psi^{*}V}\right),$$
then \eqref{eqn : gauge invariance StS}  readily follows from  $ \Psi|_{O}=\mathrm{Id}$.
\end{proof}

\section{Runge approximation}\label{SRA}
A tool to establish the Runge approximation is the following unique continuation principle of distributions for fractional elliptic equations, which is the distributional version of the entanglement principle  \cite[Theorem 1.8]{feizmohammadi2024calderonproblemfractionalschrodinger}. 
 	\begin{Lm}[Unique continuation principle]\label{UCP}
	Let $\left(M,g\right)$ be a connected closed Riemannian $n$-manifold and $O\subset M$ a nonempty connected open set. Suppose there exist $N \in \N_{+}$, $\{r_{k}\}_{k=1}^{N} \subset \R$,  $\{\alpha_{k}\}_{k=1}^{N}\subset\left(0,\infty\right)\setminus\mathbb{N}$  and  $  \{v_k\}_{k=1}^N \subset H^{-r_{k}}\left(M\right)$ such that
\begin{itemize}
	\item every $  v_{k} $ vanishes  in $O$, 
	\item every distinct pair $\alpha_l$ and $\alpha_k$ satisfy	\begin{equation*}
		\alpha_{l}-\alpha_{k}\notin \mathbb{Z}  ,
	\end{equation*}
	\item the following identity holds	\begin{align}\label{UCP_sum}
		\sum\limits_{k=1}^{N}\left(\left(-\Delta_{g}\right)^{\alpha_{k}}v_{k}\right)|_{O}&=0.
	\end{align}
\end{itemize}
Then every $v_{k}$ vanishes in $M$.
\end{Lm}
For the subsequent proof, the following lemma \cite[Proposition 3.1]{feizmohammadi2024calderonproblemfractionalschrodinger} will be particularly useful.
\begin{Lm}\label{UCP_LM}
	Let $M,g,\{\alpha_{k}\}_{k=1}^{N}$ be as in Lemma \ref{UCP}. Assume that for $\{f_{k}\}_{k=1}^{N}\subset C^{\infty}(0,\infty)$, there exist constants $c,\delta >0$ such that for all $k=1,\cdots,N$,
	\begin{equation*}
		|f_{k}(t)|\leq ce^{-\delta t}, \quad t\in [1,\infty) \quad \text{and} \quad |f_{k}(t)|\leq ce^{-\frac{\delta}{t}}, \quad t\in (0,1).
	\end{equation*}
	If there exists $l\in \N$ such that
	\begin{equation*}
		\sum\limits_{k=1}^{N}\Gamma(m+1+\alpha_{k})\int_{0}^{\infty}f_{k}(t)\frac{dt}{t^{m}} = 0, \quad \text{for all} \quad m=l,l+1,\cdots,
	\end{equation*}
then $f_{k}(t)=0, \: t\in (0,\infty)$ for all $k=1,\cdots, N$.
\end{Lm}

\begin{proof}[Proof of Lemma \ref{UCP}]
We follow the strategy of the proof of \cite[Theorem 1.8]{feizmohammadi2024calderonproblemfractionalschrodinger}. Since we deal with distributions, we have to work with the pairing of distributions with test functions, which causes additional arguments to verify the conditions in Lemma \ref{UCP_LM}.

By functional calculus, we split
\begin{equation*}
	(-\Delta_{g})^{\alpha_{k}}v_{k} = (-\Delta_{g})^{\widetilde{\alpha_{k}}}(-\Delta_{g})^{m_{k}}v_{k} 
\end{equation*}
where $\alpha_{k}=\widetilde{\alpha_{k}}+m_{k}$ with $\widetilde{\alpha_{k}}\in (0,1)$, $m_{k}\in \N$.
With unique continuation principle for Laplacian, it suffices to prove the lemma for $\{\alpha_{k}\}_{k=1}^{N}\subset (0,1)$.

We take open subsets $\omega, \omega_{1}$ of $O$ such that $\omega \Subset \omega_{1}\Subset O$ and $\phi\in C_{c}^{\infty}(O)$ with $\phi=1$ on $\omega_{1}$. For any fixed $\varphi\in C_{c}^{\infty}(\omega)$, \eqref{UCP_sum} leads to that for any $m\in \N_{+}$,
\begin{equation}\label{UCP_E0}
\sum\limits_{k=1}^{N}	\left(v_{k}, (-\Delta_{g})^{\alpha_{k}}\Delta_{g}^{m+1}\varphi   \right) = \left(\sum\limits_{k=1}^{N}(-\Delta_{g})^{\alpha_{k}}v_{k},\Delta_{g}^{m+1}\varphi\right) = 0.
\end{equation}
By functional calculus, 
\begin{equation*}
	(-\Delta_{g})^{\alpha_{k}} = \frac{1}{\Gamma(-\alpha_{k})}\int_{0}^{\infty}(e^{t\Delta_{g}}-1)\frac{dt}{t^{1+\alpha_{k}}},
\end{equation*}
it follows from \eqref{UCP_E0} and the support conditions of $\varphi$ and $v_{k}$ that
\begin{equation*}
	\sum\limits_{k=1}^{N}\left(v_{k}, \frac{1}{\Gamma(-\alpha_{k})}\int_{0}^{\infty}e^{t\Delta_{g}}\Delta_{g}^{m+1}\varphi\frac{dt}{t^{1+\alpha_{k}}}\right)=0.
\end{equation*}
We claim that
\begin{equation}\label{UCP_E1}
 \sum\limits_{k=1}^{N}\left(v_{k},\frac{1}{\Gamma(-\alpha_{k})}\int_{0}^{\infty}e^{t\Delta_{g}}\Delta_{g}^{m+1}\varphi\frac{dt}{t^{1+\alpha_{k}}}\right)= \sum\limits_{k=1}^{N}\left(v_{k}, \frac{\kappa_{k,m}}{\Gamma(-\alpha_{k})}\int_{0}^{\infty}e^{t\Delta_{g}}\Delta_{g}\varphi\frac{dt}{t^{m+1+\alpha_{k}}}\right)
\end{equation}
with $\kappa_{k,m}=\frac{\Gamma(m+1+\alpha_{k})}{\Gamma(1+\alpha_{k})}$. Equation \eqref{UCP_E1} is proved by integration by parts. We need to show that the endpoint terms vanish. First we make integration by parts once
\begin{equation}\label{UCP_E2}
	\begin{aligned}
& \left(v_{k}, \frac{1}{\Gamma(-\alpha_{k})}\int_{0}^{\infty}e^{t\Delta_{g}}\Delta_{g}^{m+1}\varphi\frac{dt}{t^{1+\alpha_{k}}}\right) \\& = \left(v_{k}, \frac{1}{\Gamma(-\alpha_{k})}\int_{0}^{\infty}\partial_{t}(e^{t\Delta_{g}}\Delta_{g}^{m}\varphi)\frac{dt}{t^{1+\alpha_{k}}}\right) \\
	&= \left(v_{k}, \frac{1}{\Gamma(-\alpha_{k})}\frac{1}{t^{1+\alpha_{k}}}e^{t\Delta_{g}}\Delta_{g}^{m}\varphi\Big|_{0}^{\infty} \right) + \left(v_{k}, \frac{1+\alpha_{k}}{\Gamma(-\alpha_{k})}\int_{0}^{\infty} e^{t\Delta_{g}}\Delta_{g}^{m}\varphi\frac{dt}{t^{2+\alpha_{k}}} \right).\\
	\end{aligned}
\end{equation}
By Sobolev embedding, we have
\begin{equation}\label{UCP_E31}
\left\|e^{t\Delta_{g}}\Delta_{g}^{m}\varphi\right\|_{L^{\infty}(M)}  \leq C\left\|e^{t\Delta_{g}}\Delta_{g}^{m}\varphi \right\|_{H^{s}(M)}, \quad \mbox{ $s>\frac{n}{2}$}.
\end{equation}
By the spectral theorem, there holds for any $q\in \R$ that
\begin{equation}\label{UCP_E3}
	\begin{aligned}
\left\|e^{t\Delta_{g}}\Delta_{g}^{m}\varphi \right\|_{H^{q}(M)} &\leq C \left\|e^{t\Delta_{g}}(1-\Delta_{g})^{\frac{q}{2}}\Delta_{g}^{m}\varphi  \right\|_{L^{2}(M)} \\
&\leq Ce^{-\lambda_{1} t} \left\|(1-\Delta_{g})^{\frac{q}{2}}\Delta_{g}^{m}\varphi \right\|_{L^{2}(M)}  \leq Ce^{-\lambda_{1} t}\left\|\Delta_{g}^{m}\varphi \right\|_{H^{q}(M)}.
\end{aligned}
\end{equation}
For $t\rightarrow \infty$, from \eqref{UCP_E31} and \eqref{UCP_E3}, 
it follows that
\begin{equation}\label{UCP_E32}
\left(v_{k}, \lim\limits_{t\rightarrow \infty}\frac{1}{\Gamma(-\alpha_{k})}\frac{1}{t^{1+\alpha_{k}}}e^{t\Delta_{g}}\Delta_{g}^{m}\varphi \right) =0.
\end{equation}
For $t\rightarrow 0$, the support conditions of $\varphi$ and $v_{k}$ imply that for any $l\in \N$,
\begin{equation*}
	\partial_{t}^{l}	\left(v_{k},e^{t\Delta_{g}}\Delta_{g}^{m}\varphi\right) = \left(v_{k} , e^{t\Delta_{g}}\Delta_{g}^{l+m}\varphi\right) \longrightarrow \left(v_{k},\lim\limits_{t\rightarrow 0}e^{t\Delta_{g}}\Delta_{g}^{l+m}\varphi\right) = \left(v_{k},\Delta_{g}^{m+l}\varphi\right) =  0,
\end{equation*}
which leads to that
\begin{equation*}
	\left(v_{k},e^{t\Delta_{g}}\Delta_{g}^{m}\varphi\right) = o\left(t^{2+\alpha_{k}}\right)  \quad  \text{as} \quad t\rightarrow 0.
\end{equation*}
Thus
\begin{equation}\label{UCP_E4}
	\left(v_{k},\lim\limits_{t\rightarrow 0}\frac{1}{\Gamma(-\alpha_{k})}\frac{1}{t^{1+\alpha_{k}}}e^{t\Delta_{g}}\Delta_{g}^{m}\varphi\right) =\frac{1}{\Gamma(-\alpha_{k})} \lim\limits_{t\rightarrow 0}\frac{1}{t^{1+\alpha_{k}}}\left(v_{k}, e^{t\Delta_{g}}\Delta_{g}^{m}\varphi\right) = 0.
\end{equation}
Combining \eqref{UCP_E32} and \eqref{UCP_E4} proves that the endpoint terms in \eqref{UCP_E2} vanish. Iterating this argument, we obtain
\begin{equation*}
	\begin{aligned}
	\left(v_{k},\frac{1}{\Gamma(-\alpha_{k})}\int_{0}^{\infty}e^{t\Delta_{g}}\Delta_{g}^{m+1}\varphi\frac{dt}{t^{1+\alpha_{k}}} \right) &= \left(v_{k}, \frac{1+\alpha_{k}}{\Gamma(-\alpha_{k})}\int_{0}^{\infty}e^{t\Delta_{g}}\Delta_{g}^{m}\varphi\frac{dt}{t^{2+\alpha_{k}}}\right) \\
&= \cdots  \\
&=\left(v_{k},\frac{\kappa_{k,m}}{\Gamma(-\alpha_{k})}\int_{0}^{\infty}e^{t\Delta_{g}}\Delta_{g}\varphi\frac{dt}{t^{m+1+\alpha_{k}}}\right)
	\end{aligned}
\end{equation*}
which is \eqref{UCP_E1}.

Next we show that the RHS of \eqref{UCP_E1} can be rewritten as
\begin{equation}\label{UCP_E5}
0= \sum\limits_{k=1}^{N}\Gamma(m+1+\alpha_{k})\int_{0}^{\infty}\frac{1}{\Gamma(1+\alpha_{k})\Gamma(-\alpha_{k})}\left(v_{k},e^{t\Delta_{g}}\Delta_{g}\varphi\right)t^{-(1+\alpha_{k})}\frac{dt}{t^{m}}.
\end{equation}
Since $v_{k}$ vanishes on $O$ for each $k=1,\cdots,N$, we have
\begin{equation*}
	\left(v_{k},\int_{0}^{\infty}e^{t\Delta_{g}}\Delta_{g}\varphi\frac{dt}{t^{m+1+\alpha_{k}}}\right)=  \left(v_{k}, (1-\phi)\int_{0}^{\infty}e^{t\Delta_{g}}\Delta_{g}\varphi \frac{dt}{t^{m+1+\alpha_{k}}}\right).
\end{equation*}
Therefore, it suffices to prove that 
\begin{equation}\label{UCP_E6}
	(1-\phi)e^{t\Delta_{g}}\Delta_{g}\varphi \frac{1}{t^{m+1+\alpha_{k}}} \in C^{\infty}\left(M,\mathcal{S}\left( \R_{+}\right)\right)
\end{equation}
and whence  Fubini's Theorem gives that
\begin{equation*}
	\begin{aligned}
	\left(v_{k},(1-\phi)\int_{0}^{\infty}e^{t\Delta_{g}}\Delta_{g}\varphi\frac{1}{t^{m+1+\alpha}}dt \right)& = \int_{0}^{\infty}\left(v_{k},(1-\phi)e^{t\Delta_{g}}\Delta_{g}\varphi\right) \frac{dt}{t^{m+1+\alpha_{k}}}  \\
	&= \int_{0}^{\infty}\left(v_{k},e^{t\Delta_{g}}\Delta_{g}\varphi\right) \frac{dt}{t^{m+1+\alpha_{k}}}
	\end{aligned}
\end{equation*}
which shows \eqref{UCP_E5}. 
We estimate the derivatives of $(1-\phi)e^{t\Delta_{g}}\Delta_{g}\varphi\frac{1}{t^{m+1+\alpha_{k}}}$ for arbitrary order $l\in \N$,
\begin{equation}\label{UCP_E50}
	\begin{aligned}
		&\left\| \partial_{t}^{l}\left((1-\phi)e^{t\Delta_{g}}\Delta_{g}\varphi\frac{1}{t^{m+1+\alpha_{k}}}\right)\right\|_{L^{\infty}(M\setminus\omega_{1})}	\\
		&=
		\left\|(1-\phi)\sum\limits_{j=0}^{l}\binom{l}{j}\partial_{t}^{j}\left(e^{t\Delta_{g}}\Delta_{g}\varphi\right)\partial_{t}^{l-j}\left(\frac{1}{t^{m+1+\alpha_{k}}} \right) \right\|_{L^{\infty}(M\setminus \omega_{1})} \\
		& \leq C\sum\limits_{j=0}^{l}\binom{l}{j} \left\| e^{t\Delta_{g}}\Delta_{g}^{j+1}\varphi \frac{(-1)^{l-j}(m+1+\alpha)\cdots(m+l-j+\alpha_{k})}{t^{m+1+l-j+\alpha_{k}}} \right\|_{L^{\infty}(M\setminus\omega_{1})} \\
		&\leq C_l \sum\limits_{j=0}^{l}	\left\|e^{t\Delta_{g}}\Delta_{g}^{j+1}\varphi  \right\|_{L^{\infty}(M\setminus \omega_{1})}\frac{1}{t^{m+1+l-j+\alpha_{k}}}.
	\end{aligned}
\end{equation}
For $t\rightarrow \infty$, we derive, from \eqref{UCP_E31}, \eqref{UCP_E3} and \eqref{UCP_E50}, that 
\begin{equation}\label{UCP_E51}
\left\| \partial_{t}^{l}\left((1-\phi)e^{t\Delta_{g}}\Delta_{g}\varphi\frac{1}{t^{m+1+\alpha_{k}}}\right)\right\|_{L^{\infty}(M\setminus\omega_{1})}	\lesssim e^{-\lambda_{1}t}.
\end{equation}
For $t\rightarrow 0$, it follows from \eqref{UCP_E50} and Lemma \ref{Li-Yau es} that there exists $\tilde{c}_{k,l}>0$,  depending on $\text{dist}_{g}(\supp(\varphi), M\setminus \omega_{1})$, $\alpha_{k}$ and $l$, such that
\begin{equation}\label{UCP_E52}
		\left\|\partial_{t}^{l}\left((1-\phi)e^{t\Delta_{g}}\Delta_{g}\varphi\frac{1}{t^{m+1+\alpha_{k}}}\right) \right\|_{L^{\infty}(M\setminus\omega_{1})}   \lesssim e^{-\frac{\tilde{c}_{k,l}}{t}}.
\end{equation}
Combining \eqref{UCP_E51} and \eqref{UCP_E52} proves \eqref{UCP_E6}, and then \eqref{UCP_E5}.

Next we denote the integrand in \eqref{UCP_E5} by 
\begin{equation*}
	f_{k}(t) :=\frac{1}{\Gamma(1+\alpha_{k})\Gamma(-\alpha_{k})} \left(v_{k},e^{t\Delta_{g}}\Delta_{g}\varphi\right)t^{-(1+\alpha_{k})}\in C^{\infty}\left((0,\infty)\right).
\end{equation*}
We prove $f_{k}(t) \equiv 0$, for all $k=1,\cdots,N$. 
For $t\rightarrow \infty$, \eqref{UCP_E3} directly gives
\begin{equation}\label{UCP_E8}
	\begin{aligned}
|f_{k}(t)| \leq C\left| \left(v_{k}, e^{t\Delta_{g}}\Delta_{g}\varphi\right)\right| t^{-(1+\alpha_{k})} \leq C\|v_{k} \|_{H^{-r_{k}}(M)} \|e^{t\Delta_{g}}\Delta_{g}\varphi \|_{H^{r_{k}}(M)}t^{-(1+\alpha_{k})}\lesssim e^{-\lambda_{1}t}.
\end{aligned}
\end{equation}
On the other hand, for the compact set $M\setminus\omega_{1} \subset M$ there exists a constant $C>0$ and integer $K\in \N_{+}$ such that for multi-index $\gamma\in \N^{n}$,
\begin{equation*}
	\left|(v_{k},v)\right|\leq C \sum\limits_{|\gamma|=0}^{K}\|\partial_{\gamma} v\|_{L^{\infty}(M\setminus\omega_{1})}, \quad \text{for all}\: v\in C_{c}^{\infty}(M)\: \text{with} \: \supp(v)\subset M\setminus\omega_{1}. 
\end{equation*} 
For $t\rightarrow 0$,  we deduce that
 \begin{equation}\label{UCP_E9}
	\begin{aligned}
	\left|\left(v_{k}, e^{t\Delta_{g}}\Delta_{g}\varphi \right)\right| &= 	\left|\left(v_{k},(1-\phi)e^{t\Delta_{g}}\Delta_{g}\varphi\right)\right| \\
&	\leq C\sum\limits_{|\gamma|=0}^{K} \left\|\partial_{\gamma}\left((1-\phi)e^{t\Delta_{g}}\Delta_{g}\varphi\right) \right\|_{L^{\infty}(M\setminus\omega_{1})} \\
		&\leq C\sum\limits_{|\gamma|=0}^{K}\left\|\partial_{\gamma}\left(e^{t\Delta_{g}}\Delta_{g}\varphi\right)\right\|_{L^{\infty}(M\setminus\omega_{1})}\\
		& \leq C\sum\limits_{|\gamma|=0}^{K} \left\|\partial_{\gamma}\left(e^{t\Delta_{g}}\Delta_{g}\varphi\right) \right\|_{H^{\lfloor \frac{n}{2}\rfloor+1}(M\setminus\omega_{1})} \\
		& \leq C\sum\limits_{|\gamma|=0}^{K}\left\| e^{t\Delta_{g}}\Delta_{g}\varphi \right\|_{H^{\lfloor\frac{n}{2}\rfloor+1+|\gamma|}(M\setminus\omega_{1})} \\
		& \leq C\sum\limits_{|\gamma|=0}^{K}\left\|e^{t\Delta_{g}}\Delta_{g}\varphi \right\|_{H^{2\lceil  {\frac{n}{4}+\frac{1}{2}+\frac{|\gamma|}{2}}  \rceil}(M\setminus\omega_{1})} \\
		& \leq C\sum\limits_{|\gamma|=0}^{K}\left( \left\|e^{t\Delta_{g}}\Delta_{g}\varphi \right\|_{L^{2}(M\setminus\omega_{1})} + \left\|\Delta_{g}^{\lceil  {\frac{n}{4}+\frac{1}{2}+\frac{|\gamma|}{2}}  \rceil}e^{t\Delta_{g}}\Delta_{g}\varphi \right\|_{L^{2}(M\setminus\omega_{1})}\right) \\
		& \leq Ce^{-\frac{c}{t}}\sum\limits_{|\gamma|=0}^{K}\left( \left\|\Delta_{g}\varphi \right\|_{L^{2}{(\omega)}}+ \left\|\Delta_{g}^{\lceil  {\frac{n}{4}+\frac{1}{2}+\frac{|\gamma|}{2}}  \rceil+1}\varphi \right\|_{L^{2}(\omega)} \right),
	\end{aligned}
\end{equation}
where the last inequality is due to Lemma \ref{LM_DG} with $c>0$ depending on $\text{dist}_{g}(\overline{\omega}, M\setminus\omega_{1})>0$.

From \eqref{UCP_E9}, it follows that 
\begin{equation}\label{UCP_E11}
	\begin{aligned}
\left|f_{k}(t)\right|\lesssim e^{-\frac{c}{t}}t^{-(1+\alpha_{k})} \lesssim e^{-\frac{c_{k}}{t}}, \quad t\in (0,1)
	\end{aligned}
\end{equation}
where $c_{k}>0$ is a constant slightly smaller than $c$.

By taking $\delta=\min\{\lambda_{1}, c_{1},\cdots, c_{N}\}$, we have proved all required conditions in Lemma \ref{UCP_LM}, i.e. \eqref{UCP_E5} for all $m\in \N_{+}$, \eqref{UCP_E8} and \eqref{UCP_E11}. Then we have  
\begin{equation*}
	f_{k}(t)\equiv0 ,\: t\in (0,\infty) \quad  \text{for all} \quad k=1,\cdots,N,
\end{equation*}
which implies 
\begin{equation*}
	\left(e^{t\Delta_{g}}\Delta_{g}v_{k},\varphi\right)=\left(v_{k},e^{t\Delta_{g}}\Delta_{g}\varphi\right) \equiv 0, \quad t>0.
\end{equation*}
By the arbitrariness of $\varphi\in C_{c}^{\infty}(\omega)$, we obtain that 
\begin{equation}
	e^{t\Delta_{g}}\Delta_{g}v_{k}(x)=0,\quad t>0, \: x\in \omega \quad \text{for all} \: k=1,\cdots,N.
\end{equation}
It follows from Lemma \ref{LUCP} that 
\begin{equation*}
		e^{t\Delta_{g}}\Delta_{g}v_{k}(x)\equiv0,\quad t>0, \: x\in M \quad \text{for all} \: k=1,\cdots,N.
\end{equation*}
Taking the limit $t\rightarrow 0$ yields
\begin{equation*}
	\Delta_{g}v_{k} = 0 \quad \text{in} \quad H^{-r_{k}-2}(M).
\end{equation*}
Since $v_{k}$ vanishes on $O$, we deduce $v_{k}=0$ in $H^{-r_{k}}(M)$ by unique continuation for elliptic equations.
\end{proof}

Now we are ready to prove the useful Runge approximation for coefficient inversion.
\begin{proof}[Proof of Theorem \ref{DR}]	
The proof uses techniques from \cite[Lemma 3.10]{feizmohammadi2024calderonproblemfractionalschrodinger}.
Recall the set $R$ defined in \eqref{DR1}. In view of the spanning criterion in \cite[Theorem 8, Page 77]{MR1892228}, we show that if $F\in\left(\widetilde{H}^{\alpha}\left(M\setminus\bar{U}\right)\right)^{*} = H^{-\alpha}(M\setminus \bar{U})$ satisfies $F\left(v\right)=0$ for all $v\in R$, then $F\equiv 0$. The proof works verbatim for $R^{*}$. To apply Lemma \ref{UCP}, we extend such $F$ to $\widetilde{F}\in \left(H^{\alpha}\left(M\right)\right)^{*}$ by 
$$ \widetilde{F}(u) = F(u\mathbbm{1}_{M\setminus U}), \quad \forall u\in H^{\alpha}(M).$$

If $\text{Ker}\left(L_{ g,b,V}\right)$ is trivial, i.e. $\text{Ker}\left(L_{ g,b,V}\right)=\{0\},$ there exists a solution $\phi\in H^{\alpha}(M)$ to \begin{align*}L_{ g,b,V}^{*}\left(\phi\right) = \widetilde{F}\end{align*}  in the sense of \eqref{SES}.
By definition, it also implies that $$H_{ g,b,V}^{U}=H_{ g,b,V}^{U,*}=C_{c}^{\infty}\left(U\right).$$ 	
For any $f \in C_{c}^{\infty}\left(U\right)$, there is a unique solution $u^f\in C^{\infty}(M)$ to \eqref{FL_E_PP} modulo $\text{Ker}(L_{g,b,V})$. Then it follows that
\begin{equation}
	0= F(u^{f}\mathbbm{1}_{M\setminus U}) =  \left(\widetilde{F},u^{f}\right)_{H^{-\alpha}(M),H^{\alpha}(M)} = \left(L_{ g,b,V}^{*}\left(\phi\right),u^{f}\right)_{L^{2}\left(M\right)}= \left(\phi,f\right)_{L^{2}\left(M\right)}. 
\end{equation}
Since $f$ is arbitrarily chosen in $C_{c}^{\infty}(U)$ and $\widetilde{F}|_{U}=0$, we see $\phi|_{U}=0$ and whence $\left(-\Delta_{g}\right)^{\alpha}\phi|_{U}=0$.  Lemma \ref{UCP} then guarantees that $\phi\equiv 0$ on $M$. This leads to that $\widetilde{F}\equiv 0 $ and $F\equiv 0$.

If $\text{Ker}\left(L_{ g,b,V}\right)$ is non-trivial, $\mathrm{Ker}(L_{g,b,V})$ is of finite dimension by Proposition \ref{WP}. Therefore, we may set
an $L^{2}\left(M\right)$-orthonormal basis $\{\zeta_{k}\}_{k=1}^{N}$  for $\text{Ker}\left(L_{ g,b,V}\right)$ with $ N\in\N_{+}$. Consider the projection onto the kernel of $L_{g,b,V}$  
\begin{align*}\pi: H^{-\alpha}\left(M\right) &\longrightarrow \text{Ker}\left(L_{ g,b,V}\right) \\
	\Xi &\longmapsto \sum\limits_{k=1}^{N} \left(\Xi,\zeta_{k}\right)_{H^{-\alpha}\left(M\right),H^{\alpha}\left(M\right)}\zeta_{k}.
\end{align*} In particular, $$\pi |_{\text{Ker}\left(L_{ g,b,V}\right)} = \mathrm{Id}.$$
Then for any $v\in \text{Ker}\left(L_{ g,b,V}\right)$ there holds that 
	\begin{align*}
		\lefteqn{\left(\widetilde{F}-\pi\left(\widetilde{F}\right), v\right)_{H^{-\alpha}\left(M\right),H^{\alpha}\left(M\right)}}  \\
		&= \left(\widetilde{F},v\right)_{H^{-\alpha}\left(M\right),H^{\alpha}\left(M\right)} -\sum\limits_{k=1}^{N} \left(\widetilde{F},\zeta_{k}\right)_{H^{-\alpha}\left(M\right),H^{\alpha}\left(M\right)} \left(\zeta_{k},v\right)_{H^{-\alpha}\left(M\right),H^{\alpha}\left(M\right)} \\
		&=\left(\widetilde{F},v-\sum\limits_{k=1}^{N}\left(v,\zeta_{k}\right)_{L^{2}\left(M\right)}\zeta_{k}\right)_{H^{-\alpha}\left(M\right),H^{\alpha}\left(M\right)} =0.
	\end{align*} 
This implies that there exists some solution $\phi\in H^{\alpha}\left(M\right)$ to 
\begin{equation}\label{eqn : L* phi = F - pi F}
	L_{ g,b,V}^{*}\left(\phi\right) = \widetilde{F}-\pi\left(\widetilde{F}\right)
\end{equation} in the sense of \eqref{SES} by Proposition \ref{WP}.

Recall, from Corollary \ref{well-defined StS}, that for any $f\in H_{ g,b,V}^{U}$ there exists a unique $u^f\in C^{\infty}(M)$ such that  
\begin{align*}L_{  g, b, V}\left(u^f\right)  = f \quad \mbox{and} \quad \left(u^{f},v\right)_{L^{2}\left(M\right)} =0,\quad \forall v\in \text{Ker}\left(L_{ g,b,V}\right).\end{align*}
It follows, from $u^{f}\mathbbm{1}_{M\setminus U}\in R$, that
\begin{equation*} 
	\begin{aligned}
		0&=F\left(u^{f}\mathbbm{1}_{M\setminus U}\right)  =   \left(\widetilde{F},u^{f}\right)_{H^{-\alpha}\left(M\right),H^{\alpha}\left(M\right)} =  \left(\widetilde{F}-\pi\left(\widetilde{F}\right),u^{f}\right)_{H^{-\alpha}\left(M\right),H^{\alpha}\left(M\right)}  
	\end{aligned}
\end{equation*}
Combining this with \eqref{eqn : L* phi = F - pi F} gives
\begin{equation}\label{REE1}
	\begin{aligned}
		0 
		&= \left(L_{ g,b,V}^{*}\left(\phi\right),u^{f}\right)_{L^{2}\left(M\right)} = \left(\phi,f\right)_{L^{2}\left(M\right)}=\left(\phi,f\right)_{L^{2}\left(U\right)}.
	\end{aligned}
\end{equation}

Since $L_{g,b,V}^{*}$ is Fredholm, $\mathrm{Ker}(L_{g,b,V}^{*})$ is also of finite dimension. Suppose $\{\eta_{k}\}_{k=1}^{N}$ is an $L^{2}\left(M\right)$-orthonormal basis for $\text{Ker}\left(L_{ g,b,V}^{*}\right)$ and write $$W:=\text{Span}\{\eta_{k}|_{U}\}_{k=1}^{N}.$$
We claim that \begin{equation}\label{OSI}
\varphi:=\phi|_{U}\in W.
\end{equation}
This is proved as in \cite[Lemma 3.10 Case II]{feizmohammadi2024calderonproblemfractionalschrodinger}.
In light of the orthogonal splitting $L^{2}\left(U\right) = W\oplus W^{\perp}$, we write $\phi|_{U} = \varphi + \phi_{0}$ with $\varphi\in W$ and $\phi_{0}\in W^{\perp}$. Then \eqref{OSI} reduces to $\phi_{0}\equiv 0.$
Let $\{h_{l}\}_{l=1}^{\infty}\subset C_{c}^{\infty}\left(U\right)$ be such that 
\begin{equation}\label{PI2}
	\left\|h_{l}-\phi_{0}\right\|_{L^{2}\left(U\right)}\rightarrow 0\quad\mbox{as $l\rightarrow \infty$}.
\end{equation}
With the fact that $\left(\phi_{0},\eta_{k}\right)_{L^{2}\left(U\right)}=0,\, \forall k=1,\cdots,N$, we deduce that
\begin{equation}\label{PI3}
	\lim\limits_{l\rightarrow \infty}\left(h_{l},\eta_{k}\right)_{L^{2}\left(U\right)} = 0,\,\forall k=1,\cdots,N.
\end{equation}
Lemma \ref{RM} gives that 
\begin{equation}\label{SSF}
	f_{l} := h_{l}- \sum\limits_{k=1}^{N}\left(h_{l},\eta_{k}\right)_{L^{2}\left(M\right)}\theta_{k}\in H_{ g,b,V}^{U}.
	\end{equation}
From \eqref{REE1} and $f_{l}\perp \text{Ker}(L_{g,b,V}^{*})$, we obtain 
\begin{equation*}\label{REE2}
	\left(f_{l},\phi_{0}\right)_{L^{2}\left(U\right)}=\left(f_{l},\phi\right)_{L^{2}\left(U\right)}-\left(f_{l},\varphi\right)_{L^{2}\left(U\right)}=0.
\end{equation*}
Together with \eqref{PI2}, \eqref{PI3} and \eqref{SSF}, we obtain $\left(\phi_{0},\phi_{0}\right)_{L^{2}\left(U\right)}=0,$ and thus \eqref{OSI} holds.

Now we extend $\varphi\in\text{Span}\{\eta_{k}|_{U}\}_{k=1}^{N}$ to $\widetilde{\psi}\in \text{Ker}\left(L_{ g,b,V}^{*}\right)$ on $M$.
Denote by $\psi:=\phi-\widetilde{\psi}\in H^{\alpha}(M)$ the solution to
\begin{equation*}
	\left\{
	\begin{aligned}
	L_{ g,b,V}^{*} \psi &= \widetilde{F}-\pi\left(\widetilde{F}\right) && \mbox{on $M$}, \\
		\psi &=0 && \mbox{on $U$}.
	\end{aligned}
	\right.
\end{equation*}
Applying $L_{ g,b,V}$ to this equation gives
\begin{equation*}
	\left\{
	\begin{aligned}
		L_{ g,b,V}\left(L_{ g,b,V}^\ast\left(\psi\right)\right) &= L_{ g,b,V}\widetilde{F} && \mbox{on $M$,} \\
		\psi &=0 && \mbox{on $U$.} 
	\end{aligned}
	\right.
\end{equation*}
Since $b$ and $V$ both vanish on $U$, we see that $\psi$ solves
\begin{equation*}
	\left(-\Delta_{g}\right)^{2\alpha}\psi +\left(-\Delta_{g}\right)^{\alpha}\nabla_{g}^{*}\left(\bar{b}\psi\right) + \left(-\Delta_{g}\right)^{\alpha}\left(\overline{V}\psi\right) = \left(-\Delta_{g}\right)^{\alpha}\widetilde{F}\quad \mbox{on $U$.}
\end{equation*}
As $\psi$ and $\widetilde{F}$ also vanish on $U$, it reduces to
\begin{align*} \left\{ \begin{aligned}
	\left(-\Delta_{g}\right)^{2\alpha}\psi+\left(-\Delta_{g}\right)^{\alpha}\left(\nabla_{g}^{*}\left(\bar{b}\psi\right)+\overline{V}\psi-\widetilde{F}\right)  &=0 && \mbox{on $U$,} \\	
	\left(\nabla_{g}^{*}\left(\bar{b}\psi\right)+\overline{V}\psi-\widetilde{F}\right) &=0 && \mbox{on $U$,} \\
		\psi &= 0 && \mbox{on $U$.} \\
	\end{aligned}\right.
\end{align*}
Since $\alpha\in \left(\frac{1}{2},1\right)$, we employ again  Lemma \ref{UCP} obtaining  that 
\begin{equation*}
	\begin{aligned}
	\nabla_{g}^{*}\left(\bar{b}\psi\right)+\overline{V}\psi-\widetilde{F} &=0\quad \text{in $H^{-\alpha}(M)$}, \\
\psi	&=0 \quad \text{in $H^{\alpha}(M)$}.
	\end{aligned}
\end{equation*}
Therefore, we get $\widetilde{F}\equiv 0$ and $F\equiv 0.$
\end{proof}

\section{Recovery of the lower order terms}\label{S4}
Before proving Theorem \ref{MR1}, we first draw the link between the local source-to-solution map and the lower order terms.
\begin{Lm}\label{AI}
	Let $M$ be a connected closed Riemannian $n$-manifold and $O\subset M$ a nonempty open set. Let $f\in H_{ g,b_{1},V_{1}}^{O}\cap H_{ g,b_{2},V_{2}}^{O}$ and $h\in H_{ g,b_{2},V_{2}}^{O,*}$. There holds
	\begin{align*}
		\lefteqn{	\left(\Lambda_{ g,b_{1},V_{1}}^{O}\left(f\right)-\Lambda_{ g,b_{2},V_{2}}^{O}\left(f\right),h \right)_{L^{2}\left(M\right)}}\\
	&= \left(\left(b_{2}-b_{1}\right)\left(\nabla_{g} u_{1}^{f}\right), u_{2}^{h,*}\right)_{L^{2}\left(M\right)} + \left(\left(V_{2}-V_{1}\right)u_{1}^{f},u_{2}^{h,*}\right)_{L^{2}\left(M\right)}.
	\end{align*}    			
\end{Lm}
\begin{proof} This lemma follows from a direct computation.
	\begin{align*}
		\lefteqn{	 \left(\Lambda_{ g,b_{1},V_{1}}^{O}\left(f\right)-\Lambda_{ g,b_{2},V_{2}}^{O}\left(f\right), h\right)_{L^{2}\left(M\right)} } \\
&= \left(u_{1}^{f}-u_{2}^{f},L_{ g,b_{2},V_{2}}^{*}(u_{2}^{h,*})\right)_{L^{2}\left(M\right)} \\
			&= \left(L_{ g,b_{2},V_{2}}(u_{1}^{f}),u_{2}^{h,*}\right)_{L^{2}\left(M\right)} - \left(L_{ g,b_{2},V_{2
			}}(u_{2}^{f}),u_{2}^{h,*}\right)_{L^{2}\left(M\right)} \\
			&= \left(L_{ g,b_{2},V_{2}}(u_{1}^{f}),u_{2}^{h,*}\right)_{L^{2}\left(M\right)} - \left(L_{ g,b_{1},V_{1}}(u_{1}^{f}),u_{2}^{h,*}\right)_{L^{2}\left(M\right)} \\
			&= \left(\left(b_{2}-b_{1}\right)(\nabla_{g} u_{1}^{f}),u_{2}^{h,*}\right)_{L^{2}\left(M\right)} + \left(\left(V_{2}-V_{1}\right)u_{1}^{f},u_{2}^{h,*}\right)_{L^{2}\left(M\right)}.
	\end{align*}       		
\end{proof}

\begin{proof}[Proof of Theorem \ref{MR1}]     		       		
 Take $U\subset O$ a nonempty open set with Lipschitz boundary. For simplicity, we denote by $u^{f}$ the unique solution to 
\begin{equation}\label{E1}
	L_{ g,b_{1},V_{1}}\left(u^{f}\right)=f, \quad f\in H_{g,b_{1},V_{1}}^{U}
\end{equation} subject to $\left(u^{f},v\right)_{L^{2}\left(M\right)}=0,\, \forall v\in \mathrm{Ker}\left(L_{g,b_{1},V_{1}}\right)$, and by $u^{h,*}$ the unique solution to 
\begin{equation}\label{E2}
	L_{ g,b_{2},V_{2}}^{*}\left(u\right)=h, \quad h\in H_{g,b_{2},V_{2}}^{U,*}
\end{equation}
subject to $\left(u^{h,*},v\right)_{L^{2}\left(M\right)}=0,\, \forall v\in \mathrm{Ker}\left(L_{g,b_{2},V_{2}}^{*}\right)$.

Applying Lemma \ref{AI}, we obtain 
	\begin{equation}\label{RPRI}
		\begin{aligned}
		0&=\left(\Lambda_{ g,b_{1},V_{1}}^{U}\left(f\right)-\Lambda_{ g,b_{2},V_{2}}^{U}\left(f\right),h\right)_{L^{2}\left(M\right)}\\
		&=\left(\left(b_{2}-b_{1}\right)(\nabla_{g} u^{f}), u^{h,*}\right)_{L^{2}\left(M\right)} +
		\left(\left(V_{2}-V_{1}\right)u^{f},u^{h,*}\right)_{L^{2}\left(M\right)}. 
		\end{aligned}
	\end{equation}	
	
	For any $\chi^\dagger\in C_{c}^{\infty}\left(M \setminus\bar{U}\right)$, we select $\chi^\ddagger\in C_{c}^{\infty}\left(M\setminus\bar{U}\right)$ such that $\chi^\ddagger=1$ on $\text{supp}\left(\chi^\dagger\right)$. By Theorem \ref{DR}, there exist sequences of solutions $\{u^{f_{k}}\}$ to \eqref{E1} and $\{u^{h_{k},*}\}$ to \eqref{E2} such that \begin{equation}\label{CS}
		\left\{
		\begin{aligned} r_{k}  := u^{f_{k}}\mathbbm{1}_{M\setminus U} - \chi^{\ddagger} &\longrightarrow 0 && \mbox{in $\widetilde{H}^{\alpha}\left(M\setminus \bar{U}\right)$}, \quad \mbox{as $k\rightarrow \infty$};\\ 
			r_{k}^{*}  := u^{h_{k},*}\mathbbm{1}_{M\setminus U} - \chi^\dagger &\longrightarrow 0 && \mbox{in $\widetilde{H}^{\alpha}\left(M\setminus\bar{U}\right)$}, \quad \mbox{as $k\rightarrow \infty$} .  \end{aligned}
			\right.
	\end{equation}
	
Through such constructions, we readily have 
	\begin{equation*} 
	\begin{aligned}
		\lefteqn{  \left| \left(\left(b_{2}- b_{1}\right)\left(\nabla_{g} u^{f_{k}}\right), u^{h_{k},*} \right)_{L^{2}\left(M\right)} \right|}\\
		&= \left|\left((b_{2}-b_{1})\left(\nabla_{g} \left(\chi^\ddagger+r_{k}\right)\right), \chi^\dagger+r_{k}^{*} \right)_{L^{2}(M)} \right| \\
		& \leq \left|\left((b_{2}-b_{2})(\nabla_{g} \chi^\ddagger),r_{k}^{*}\right)_{L^{2}\left(M\right)}\right| + \left|\left(\nabla_{g} r_{k},\chi^\dagger\left(\bar{b}_{2}-\bar{b}_{1}\right)\right)_{L^{2}\left(M\right)}\right| \\
		& \quad +  
		\left|\left(\nabla_{g} r_{k},r_{k}^{*}\left(\bar{b}_{2}-\bar{b}_{1}\right)\right)_{L^{2}\left(M\right)}\right|
	\end{aligned}
\end{equation*}
where the last inequality is due to $\chi^{\ddagger}=1$ on $\text{supp}(\chi^{\dagger})$.

Since $b_{1}, \, b_{2}$ are smooth on $M$, $b_{1}=b_{2}$ on $U$ and $0<1-\alpha<\alpha<1$, we deduce
$$
		\begin{aligned}
			\lefteqn{  \left| \left(\left(b_{2}- b_{1}\right)\left(\nabla_{g} u^{f_{k}}\right), u^{h_{k},*} \right)_{L^{2}\left(M\right)} \right|} \\
		&	\leq \left\|(b_{2}-b_{1})(\nabla_{g}\chi^{\ddagger}) \right\|_{H^{-\alpha}(M\setminus\bar{U})} \|r_{k}^{*} \|_{\widetilde{H}^{\alpha}(M\setminus\bar{U})} \\
			& \quad + \left\|\nabla_{g} r_{k} \right\|_{H^{\alpha-1}\left(M\setminus\bar{U}\right)}  \left\|\chi^\dagger\left(\bar{b}_{2}-\bar{b}_{1}\right) \right\|_{\widetilde{H}^{1-\alpha}\left(M\setminus\bar{U}\right)} \\
			& \quad + \|\nabla_{g} r_{k} \|_{H^{\alpha-1}\left(M\setminus\bar{U}\right)} \left\|r_{k}^{*}\left(\bar{b}_{2}-\bar{b}_{1}\right) \right\|_{\widetilde{H}^{1-\alpha}\left(M\setminus\bar{U}\right)} \\
			&\leq \left\|(b_{2}-b_{1})(\nabla_{g}\chi^{\ddagger})\right\|_{H^{-\alpha}(M\setminus\bar{U})} \left\|r_{k}^{*} \right\|_{\widetilde{H}^{\alpha}(M\setminus\bar{U})}  \\
			& \quad +\sqrt{n}\|r_{k} \|_{H^{\alpha}\left(M\setminus\bar{U}\right)}\left(\|b_{2}-b_{1} \|_{L^{\infty}\left(M\right)}+\left\|\left(1-\Delta_{g}\right)^{\frac{1-\alpha}{2}}\left(b_{2}-b_{1}\right) \right\|_{L^{\infty}\left(M\right)}\right)\left\|\chi^\dagger \right\|_{\widetilde{H}^{1-\alpha}\left(M\setminus\bar{U}\right)}\\
			&\quad +\sqrt{n} \|r_{k} \|_{H^{\alpha}\left(M\setminus\bar{U}\right)}  \left(\|b_{2}-b_{1}\|_{L^{\infty}\left(M\right)}+  \left\|\left(1-\Delta_{g}\right)^{\frac{1-\alpha}{2}}\left(b_{2}-b_{1}\right) \right\|_{L^{\infty}\left(M\right)}   \right)\left\|r_{k}^{*} \right\|_{\widetilde{H}^{1-\alpha}\left(M\setminus\bar{U}\right)} \\
			&\lesssim  \|r_{k}^{*} \|_{\widetilde{H}^{\alpha}(M\setminus\bar{U})} \left\|(b_{2}-b_{1})(\nabla_{g}\chi^{\ddagger})\right\|_{H^{-\alpha}(M\setminus\bar{U})}\\
			&\quad + C_{b_{1}, b_{2}}\left\|r_{k} \right\|_{\widetilde{H}^{\alpha}\left(M\setminus\bar{U}\right)}\left( \left\| \chi^\dagger \right\|_{\widetilde{H}^{\alpha}\left(M\setminus\bar{U}\right)}+ \left\|r_{k}^{*} \right\|_{\widetilde{H}^{\alpha}\left(M\setminus\bar{U}\right)} \right)
		\end{aligned}
$$
where $$
	C_{b_{1}, b_{2}} := \|b_{2}-b_{1}\|_{L^{\infty}\left(M\right)}+  \left\|\left(1-\Delta_{g}\right)^{\frac{1-\alpha}{2}}\left(b_{2}-b_{1}\right) \right\|_{L^{\infty}\left(M\right)}< \infty. 
$$
Together  with \eqref{CS}, we obtain that 
	\begin{equation}\label{PRE1} \left| \left(\left(b_{2}- b_{1}\right)\left(\nabla_{g} u^{f_{k}}\right), u^{h_{k},*} \right)_{L^{2}\left(M\right)} \right|  \longrightarrow 0, \quad \mbox{as $k\rightarrow \infty$.}	\end{equation}  
	
	On the other hand, it is easy to compute that
	\begin{equation*} 
		\begin{aligned}
		\lefteqn{	\left(\left(V_{2}-V_{1}\right) u^{f_{k}}, u^{h_{k},*}\right)_{L^{2}\left(M\right)}}\\
			&= \left(\left(V_{2}-V_{1}\right)\left(\chi^\ddagger+r_{k}\right),\chi^\dagger+r_{k}^{*}\right)_{L^{2}\left(M\setminus\bar{U}\right)}\\
			&= \left(\left(V_{2}-V_{1}\right)\chi^\ddagger,\chi^\dagger\right)_{L^{2}\left(M\setminus\bar{U}\right)} +\left(\left(V_{2}-V_{1}\right)r_{k},\chi^\dagger\right)_{L^{2}\left(M\setminus\bar{U}\right)}\\
			&+\left(\left(V_{2}-V_{1}\right)\chi^\ddagger,r_{k}^{*}\right)_{L^{2}\left(M\setminus\bar{U}\right)}+\left(\left(V_{2}-V_{1}\right)r_{k},r_{k}^{*}\right)_{L^{2}\left(M\setminus\bar{U}\right)}.
		\end{aligned}
	\end{equation*} By \eqref{CS}, it is clear that
	\begin{equation}\label{PRE2}
	\begin{aligned}
	 	\left(\left(V_{2}-V_{1}\right) u^{f_{k}}, u^{h_{k},*}\right)_{L^{2}\left(M\right)}  \longrightarrow \int_{M\setminus\bar{U}}\left(V_{2}-V_{1}\right)\overline{\chi^\dagger} dV_{g},\quad \mbox{as $k\rightarrow\infty$}.
	\end{aligned}
\end{equation}

 In summary, identity \eqref{RPRI}, together with \eqref{PRE1} and \eqref{PRE2} for $f_k\in H_{g,b_{1},V_{1}}^{U}$ and $h_k\in H_{g,b_{2},V_{2}}^{U,*}$, shows that
	\begin{equation*}
	 \int_{M\setminus\bar{U}}\left(V_{2}-V_{1}\right)\overline{\chi^\dagger}dV_{g} = 0, \quad \forall \chi^\dagger\in C_{c}^{\infty}\left(M\setminus\bar{U}\right).
	\end{equation*}
Therefore, $V_{2}$ and $V_{1}$ must agree on $M\setminus\bar{U}$ and also agree on $M\setminus\bar{O}$.

Furthermore, this in turn leads to
	\begin{equation}\label{RPRI2}
		\left(\left(b_{2}-b_{1}\right)\left(\nabla_{g} u^{f}\right),u^{h,*}\right)_{L^{2}\left(M\right)} = 0
	\end{equation}

	For any point $x_{0}\in M\setminus\bar{U}$, we now take $\chi_0\in C_{c}^{\infty}\left(W\right)\subset C_{c}^{\infty}\left(M\setminus\bar{U}\right)$ to be supported within a small neighbourhood $W$ of $x_{0}$. We adopt local coordinates $x = \left(x^1, \cdots, x^n\right)$ in $W$ and select $\chi_{0l}\in C_{c}^{\infty}\left(M\setminus\bar{U}\right)$ such that $\chi_{0l}\left(x\right)=x^{l}$ on $\text{supp}\left(\chi_0\right)\subset M\setminus\bar{U}$. By Theorem \ref{DR}, there exist sequences of solutions $\{u^{f_{0k}}\}$ and $\{u^{h_{0k},*}\}$ such that \begin{equation}\label{CS2}
				\left\{	\begin{aligned} 
			r_{0k}:= u^{f_{0k}}\mathbbm{1}_{M\setminus U} -\chi_{0l} &\longrightarrow 0 && \mbox{in $\widetilde{H}^{\alpha}\left(M\setminus \bar{U}\right)$}, \quad \mbox{as $k\rightarrow\infty$};\\
			r_{0k}^{*}:= u^{h_{0k},*}\mathbbm{1}_{M\setminus U} - \chi_{0} & \longrightarrow 0 && \mbox{in $\widetilde{H}^{\alpha}\left(M\setminus\bar{U}\right)$}, \quad \mbox{as $k\rightarrow\infty$}.
			\end{aligned}
			\right.
			\end{equation}
	
	Inserting functions $\{u^{f_{0k}}\}\,\text{and}\,\{u^{h_{0k},*}\}$ into \eqref{RPRI2}, we obtain	
	\begin{align*}0&= \left(\left(b_{2}-b_{1}\right)\left(\nabla_{g} u^{f_{0k}}\right), u^{h_{0k},*}\right)_{L^{2}\left(M\right)}\\
			 &= \left(\left(b_{2}-b_{1}\right)\left( \nabla_{g} \left(\chi_{0l}+r_{0k}\right)\right),\chi_0+r_{0k}^{*}\right)_{L^{2}\left(M\setminus\bar{U}\right)}\\
			&= \left(\left(b_{2}-b_{1}\right)\left(\nabla_{g}\chi_{0l}\right),\chi_{0}\right)_{L^{2}\left(M\setminus\bar{U}\right)} + \left(\left(b_{2}-b_{1}\right)\left(\nabla_{g} r_{0k}\right),\chi_{0}\right)_{L^{2}\left(M\setminus\bar{U}\right)} \\
			&+ \left(\left(b_{2}-b_{1}\right)\left(\nabla_{g}\chi_{0l}\right),r_{0k}^{*}\right)_{L^{2}\left(M\setminus\bar{U}\right)}
			+ \left(\left(b_{2}-b_{1}\right)\left(\nabla_{g} r_{0k}\right),r_{0k}^{*}\right)_{L^{2}\left(M\setminus\bar{U}\right)}.
	\end{align*}
	As $k\rightarrow \infty$, \eqref{CS2} gives that
	$$0=	\left(\left(b_{2}-b_{1}\right)\left(\nabla_{g} u^{f_{0k}}\right), u^{h_{0k},*}\right)_{L^{2}\left(M\right)} \longrightarrow \int_{M\setminus\bar{U}}\left(\left(b_{2}-b_{1}\right)_{k}g^{lk}\right)\overline{\chi}_{0}dV_{g}.$$	
	Since $\chi_0$ is arbitrarily chosen around $x_0$ and $l$ is arbitrary in $\{1,\cdots,n\}$, it follows that $b_{2}\left(x_{0}\right)= b_{1}\left(x_{0}\right)$. As $x_0$ is also arbitrarily chosen in $M\setminus\bar{U}$,  we conclude that  $b_{1}$ and $b_{2}$ coincide on $M\setminus\bar{U}$ and thus $b_{1}=b_{2}$ on $M\setminus \bar{O}$.	
\end{proof}

\section{Recovery of the Riemannian metric}\label{S3}
To conclude the proof of Theorem \ref{MR}, it remains to recover metric $g$. This step can be carried out by virtue of the proof in \cite[Section 4]{feizmohammadi2024calderonproblemfractionalschrodinger} (originally for fractional Schrödinger equations without drift terms, but adaptable to our framework with minimal adjustments). For the reader's convenience, we sketch the key steps of this adapted strategy below.

For $j= 1,2$, we define the solution sets $K^{(j)}$ as follows:
\begin{equation}\label{DSS}
	K^{(j)} :=\left\{u\in C^{\infty}(M): \mbox{$L_{g_{j},b_{j},V_{j}}(u)\in H_{g_{j},b_{j},V_{j}}^{O}$ and $u\perp \text{Ker}(L_{g_{j},b_{j},V_{j}})$ in $L^{2}(M)$}\right\},
\end{equation}
and denote by $u_{j}^{f}$ an element in $K^{(j)}$ such that  $L_{g_{j},b_{j},V_{j}}(u_{j}^{f})=f\in H_{g_{j},b_{j},V_{j}}^{O}$.

\subsection{From the source-to-solution map to the heat semigroup}\label{S31} 
	\begin{Prop} \label{P1} Let $\alpha, M, O, g_j, b_j, V_j, \Lambda^O_{g_j, b_j, V_j}$ be as in Theorem \ref{MR}. 
	 If $\Lambda_{ g_{1},b_{1},V_{1}}^{O}=\Lambda_{ g_{2},b_{2},V_{2}}^{O}$, then for any $f\in H_{g_{1},b_{1},V_{1}}^{O}=H_{g_{2},b_{2},V_{2}}^{O}$, there holds
		\begin{equation}\label{HSI}
			e^{t\Delta_{g_{1}}}\left(-\Delta_{g_{1}}\right)^{\alpha}u_{1}^{f}|_{O}= e^{t\Delta_{g_{2}}}\left(-\Delta_{g_{2}}\right)^{\alpha}u_{2}^{f}|_{O},\quad t>0.
		\end{equation}
	\end{Prop}
	\begin{proof}This proof is the same as in  \cite[Lemma 4.3]{feizmohammadi2024calderonproblemfractionalschrodinger}. For brevity, we use $u_{j}$ to abbreviate $u_{j}^{f}$.
	 For any integer $m>0$, the identity $$\left(\Delta_{g_{j}}^{m}\left(f-b_{j}(\nabla_{g_{j}} u_{j})-V_{j}u_{j}\right),1\right)_{L^{2}\left(M\right)}=\left(f-b_{j}(\nabla_{g_{j}} u_{j})-V_{j}u_{j},\Delta_{g_{j}}^{m}1\right)_{L^{2}\left(M\right)}=0,$$ 
	together with \eqref{IFL}, implies that $$\Delta_{g_{j}}^{m}\left(f-b_{j}(\nabla_{g_{j}}u_{j})-V_{j}u_{j}\right)\in D\left( (-\Delta_{g})^{-\alpha}\right).$$
Applying $\Delta_{g_{j}}^{m}$ to the equation
		\begin{equation*}
			\left(-\Delta_{g_{j}}\right)^{\alpha}u_{j} = f-b_{j} \left( \nabla_{g_{j}} u_{j}\right)-V_{j}u_{j}, \quad f\in H_{g_{j},b_{j},V_{j}}^{O},
		\end{equation*} gives
	\begin{equation*}
	\Delta_{g_{j}}^{m}u_{j} = \left(-\Delta_{g_{j}}\right)^{-\alpha}\Delta_{g_{j}}^{m}\left(f-b_{j}(\nabla_{g_{j}} u_{j})-V_{j}u_{j}\right).
\end{equation*}
Then it follows, from the coincidence of $\Lambda_{ g_{j},b_{j},V_{j}}^{O}$ and $H_{ g_{j},b_{j},V_{j}}^{O}$, that
		 \begin{equation}\label{TH3I}
		 	\left(-\Delta_{g_{1}}\right)^{-\alpha}\Delta_{g_{1}}^{m} \left(f-b_{1}\left(\nabla_{g_{1}} u_{1}\right)-V_{1}u_{1}\right)|_{O} = \left(-\Delta_{g_{2}}\right)^{-\alpha}\Delta_{g_{2}}^{m} \left(f-b_{2}\left(\nabla_{g_{2}} u_{2}\right)-V_{2}u_{2}\right)|_{O}.
		 \end{equation} 
By functional calculus
		 \begin{equation*}\label{TH3}
		 	\left(-\Delta_{g_{j}}\right)^{-\alpha} = \frac{1}{\Gamma\left(\alpha\right)}\int_{0}^{\infty}e^{t\Delta_{g_{j}}}\frac{dt}{t^{1-\alpha}},
		 \end{equation*}
\eqref{TH3I} leads to that for any $x\in O$,
		\begin{equation*}
			\int_{0}^{\infty}\partial_{t}^{m}\left(e^{t\Delta_{g_{1}}}\left(f-b_{1}(\nabla_{g_{1}} u_{1})-V_{1}u_{1}\right)-e^{t\Delta_{g_{2}}}\left(f-b_{2}(\nabla_{g_{2}} u_{2})-V_{2}u_{2}\right)\right)(x)\frac{dt}{t^{1-\alpha}} =0.
			\end{equation*}
		 Via a straight adaptation of the integration by parts argument in \cite[Theorem 1.1]{feizmohammadi2021fractionalanisotropiccalderonproblem}, we obtain
	\begin{equation}\label{eqn : IBP}
		\int_{0}^{\infty}\left(e^{t\Delta_{g_{1}}}\left(f-b_{1}(\nabla_{g_{1}} u_{1})-V_{1}u_{1}\right)-e^{t\Delta_{g_{2}}}\left(f-b_{2}(\nabla_{g_{2}} u_{2})-V_{2}u_{2}\right)\right)(x)\frac{dt}{t^{1+m-\alpha}} = 0.
	\end{equation}
	Setting $s:=t^{-1}$, we rewrite \eqref{eqn : IBP} as \begin{equation}\label{D_FT}
		\int_{0}^{\infty}Q\left(s;x\right)s^{m-1}ds =0,
	\end{equation}
	where $$Q(s;x) :=s^{-\alpha} \left(e^{\frac{1}{s}\Delta_{g_{1}}}\left(f-b_{1}(\nabla_{g_{1}} u_{1})-V_{1}u_{1}\right)-e^{\frac{1}{s}\Delta_{g_{2}}}\left(f-b_{2}(\nabla_{g_{2}} u_{2})-V_{2}u_{2}\right)\right)(x).$$		
	
	Now we  take an open set $\omega \Subset O$ such that $\supp \left(f\right) \cap \overline{\omega} = \emptyset.$ Since $b_{j}$ and $V_{j}$ both vanish on $O$, $$ \overline{\omega}\cap \text{supp}\left(f-b_{j}(\nabla u_{j})-V_{j}u_{j}\right)=\emptyset.$$
 By Lemma \ref{DG_E}, there holds that
			\begin{equation*}
				\|Q\left(s;\cdot \right) \|_{L^{1}\left(\omega\right)}\leq C\frac{e^{-cs}}{s^{\alpha}},\, s>0.
			\end{equation*}
			 It follows that the Fourier transform of $\mathbbm{1}_{\left(0,\infty\right)}Q$ in $s$-variable, 
			\begin{equation}\label{FT}
				\mathcal{F}\left(\mathbbm{1}_{\left(0,\infty\right)}Q\right)\left(\xi;x\right)=\int_{0}^{\infty}e^{-is\xi}Q\left(s;x\right)ds,\quad x\in \omega
				\end{equation} is holomorphic for $\Im \xi < c$. 
			In view of \eqref{D_FT},  $\mathcal{F}\left(\mathbbm{1}_{\left(0,\infty\right)}Q\right)\left(\xi;x\right)$, as well as its derivatives in $\xi$ of any order, vanishes at $\xi= 0$. The analyticity of \eqref{FT} gives  $$\mathcal{F}\left(\mathbbm{1}_{\left(0,\infty\right)}Q\right)\left(\xi;x\right)=0, \quad \Im \xi<c, \: x\in \omega.$$			
			Inverting the Fourier transform, we obtain that 
			\begin{equation*}
				Q\left(s;x\right)=0, \, s\in \left(0,\infty\right),\, x\in \omega.
			\end{equation*}			
			Hence we get
			\begin{equation*}
				e^{t\Delta_{g_{1}}}\left(f-b_{1}(\nabla_{g_{1}} u_{1})-V_{1}u_{1}\right)|_{\omega} = e^{t\Delta_{g_{2}}}\left(f-b_{2}(\nabla_{g_{2}} u_{2})-V_{2}u_{2}\right)|_{\omega},\: t>0.
				\end{equation*}
			By Lemma \ref{LUCP} and agreement of $g_{1}$ and $g_{2}$ on $O$, we get
						\begin{equation*}
				e^{t\Delta_{g_{1}}}\left(f-b_{1}(\nabla_{g_{1}} u_{1})-V_{1}u_{1}\right)|_{O} = e^{t\Delta_{g_{2}}}\left(f-b_{2}(\nabla_{g_{2}} u_{2})-V_{2}u_{2}\right)|_{O},\: t>0,
			\end{equation*}
			which is equivalent to \eqref{HSI}.
		\end{proof}
				
\subsection{Recovery of the spectral data}\label{S32}
	    Now we prove that the local heat semigroup data \eqref{RM1I} with $H_{g,0,V}^{O}$ replaced by $H_{g,b,V}^{O}$ recovers the spectral data \eqref{SPD}.
		\begin{Prop}\label{P2}
		Let $\alpha,M,O,g_{j},b_{j},V_{j}$ be as in Theorem \ref{MR}. For any $f\in H_{g_{1},b_{1},V_{1}}^{O}=H_{g_{2},b_{2},V_{2}}^{O}$, the heat semigroup data \eqref{HSI} determines
\begin{equation}\label{SDI}
	 \left\{\lambda_{k}^{\left(1\right)},\pi_{k}^{\left(1\right)}u_{1}^{f}|_{O}\right\}_{k=0}^{\infty}=\left\{\lambda_{k}^{\left(2\right)},\pi_{k}^{\left(2\right)}u_{2}^{f}|_{O}\right\}_{k=0}^{\infty},
\end{equation}
				where $\{\lambda_{k}^{\left(j\right)}\}_{k=0}^{\infty}$ is a strictly monotone sequence of eigenvalues of $-\Delta_{g_{j}}$, and $\pi_{k}^{\left(j\right)}$ is the projection onto the corresponding eigenspace of $\lambda_{k}^{(j)}$ defined in \eqref{DP}.
			\end{Prop}
       \begin{proof}This proof is the same as in \cite[Proposition 4.1]{feizmohammadi2024calderonproblemfractionalschrodinger}.
       	We rewrite \eqref{HSI} in terms of the spectral representation, 
       	\begin{equation}\label{SR}
       		\sum\limits_{k=0}^{\infty}\left(\lambda_{k}^{\left(1\right)}\right)^{\alpha}e^{-t\lambda_{k}^{\left(1\right)}}\pi_{k}^{\left(1\right)}u_{1}^{f}(x) = 	\sum\limits_{k=0}^{\infty}\left(\lambda_{k}^{\left(2\right)}\right)^{\alpha}e^{-t\lambda_{k}^{\left(2\right)}}\pi_{k}^{\left(2\right)}u_{2}^{f}(x),\quad x\in O, \, t>0.
       		\end{equation}
       		\cite[Proposition 4.1]{feizmohammadi2024calderonproblemfractionalschrodinger} proves that the series in \eqref{SR} converges uniformly in $x\in M$.
       		By Laplace transforming \eqref{SR} in $t$, we obtain
       		\begin{equation}\label{SRSS}
       			\sum\limits_{k=1}^{\infty}\frac{\left(\lambda_{k}^{\left(1\right)}\right)^{\alpha}\pi_{k}^{\left(1\right)}u_{1}^{f}\left(x\right)}{\lambda_{k}^{\left(1\right)}+z} = 	\sum\limits_{k=1}^{\infty}\frac{\left(\lambda_{k}^{\left(2\right)}\right)^{\alpha}\pi_{k}^{\left(2\right)}u_{2}^{f}\left(x\right)}{\lambda_{k}^{\left(2\right)}+z}, \quad x\in O, \, \Re z>0.
       			\end{equation}
       		Denote $\Omega_{j} := \mathbb{C}\setminus \bigcup_{k=1}^{\infty}\{-\lambda_{k}^{\left(j\right)}\}$. The series 	$$R_{j}\left(z;x\right):=\sum\limits_{k=1}^{\infty}\frac{\left(\lambda_{k}^{\left(j\right)}\right)^{\alpha}\pi_{k}^{\left(j\right)}u_{j}^{f}\left(x\right)}{\lambda_{k}^{\left(j\right)}+z}$$ is holomorphic in $\Omega_{j}$ as a function of $z$.
       		
       		Without loss of generality, assume that $\lambda_{1}^{\left(1\right)}\leq \lambda_{1}^{\left(2\right)}$. From \eqref{SRSS}, we have for any $x\in O$
       		\begin{equation}\label{REV}
       			\begin{aligned}
       				\left(\lambda_{1}^{\left(1\right)}\right)^{\alpha}\pi_{1}^{\left(1\right)}u_{1}^{f}\left(x\right) &=\lim\limits_{z\rightarrow -\lambda_{1}^{\left(1\right)}}\left(\lambda_{1}^{\left(1\right)}+z\right)R_{1}\left(z;x\right) \\
       				&=\lim\limits_{z\rightarrow -\lambda_{1}^{\left(1\right)}}\left(\lambda_{1}^{\left(1\right)}+z\right)R_{2}\left(z;x\right) \\
       				&= \left\{
       				\begin{array}{ll}
       					 0, &\lambda_{1}^{\left(1\right)}< \lambda_{1}^{\left(2\right)}; \\
       					 \left(\lambda_{1}^{\left(2\right)}\right)^{\alpha}\pi_{1}^{\left(2\right)}u_{2}^{f}\left(x\right), &\lambda_{1}^{\left(1\right)}=\lambda_{1}^{\left(2\right)}.
       					\end{array}
       				\right.
       				\end{aligned}
       			\end{equation}
       		By Lemma~\ref{SK}, there exists some $f\in H_{ g_{1},b_{1},V_{1}}^{O}$ such that $$\pi_{1}^{\left(1\right)}u_{1}^{f} = \sum\limits_{l=1}^{d_{1}^{(1)}}(u_{1}^{f},\phi_{1,l}^{(1)})\phi_{k,1}^{(1)}\neq 0\quad \mbox{on $O$}.$$ Therefore, \eqref{REV} implies that $$\lambda_{1}^{\left(1\right)}=\lambda_{1}^{\left(2\right)}\quad \text{and} \quad  			\pi_{1}^{\left(1\right)}u_{1}^{f}|_{O}=\pi_{1}^{\left(2\right)}u_{2}^{f}|_{O}.$$  In the end, \eqref{SDI} is proved by induction.
       	   	\end{proof}

       	   	\begin{Prop}\label{P3}
       	   	Let $\alpha, M, O, g_{j}, b_{j}, V_{j}$ be as in Theorem \ref{MR}. The information \eqref{SDI} recovers the spectral data consisting of eigenfunctions $\left\{\{\phi_{k,l}^{(j)}\}_{l=1}^{d_{k}^{(j)}}\right\}_{k=0}^{\infty}$ with eigenvalues $\{\lambda_k^{(j)}\}_{k=0}^{\infty}$ of $-\Delta_{g_{j}}:$
       	   	\begin{equation}\label{SDI2}
       	   		 \left\{\lambda_{k}^{\left(1\right)},\{\phi_{k,l}^{\left(1\right)}|_{O}\}_{l=1}^{d_{k}^{(1)}}\right\}_{k=0}^{\infty} = \left\{\lambda_{k}^{\left(2\right)},\{\phi_{k,l}^{\left(2\right)}|_{O}\}_{l=1}^{d_{k}^{(2)}}\right\}_{k=0}^{\infty}.
       	   	\end{equation}
       	   		In addition, $\left\{\{\phi_{k,l}^{\left(j\right)}\}_{l=1}^{d_{k}^{\left(j\right)}}\right\}_{k=0}^{\infty}$ for $j=1,2$ can be made an $L^{2}\left(M\right)$-orthonormal basis simultaneously.
       	   		\end{Prop}       	   		
       	   	\begin{proof}This proof is the same as in  \cite[Theorem 1.11 and Proposition 4.1]{feizmohammadi2024calderonproblemfractionalschrodinger}.
       	   		We first prove \eqref{SDI2}.
       	   		 For any $k\in\mathbb{N}$ and $j=1,2$, denote $$S_{k}^{(j)} := \mathrm{Span}\left\{\pi_{k}^{\left(j\right)}u_{j}^{f}: f\in H_{ g_{j},b_{j},V_{j}}^{O}\right\} \subset \text{Ker}\left(-\Delta_{g_{j}}-\lambda_{k}^{\left(j\right)}\right).$$
       	   		We claim that $S_{k}^{(j)}$ is just the eigenspace of $-\Delta_{g_{j}}$ associated with $\lambda_{k}^{(j)}$. That is,
       	   		\begin{equation}\label{DR3}
       	   			S_{k}^{(j)}= \mathrm{Ker} \left(-\Delta_{g_{j}}-\lambda_{k}^{\left(j\right)}\right).
       	   		\end{equation}   	   		
       	   		If $\text{Ker}\left(-\Delta_{g_{j}}-\lambda_{k}^{\left(j\right)}\right)\setminus S_{k}^{(j)}\neq \emptyset$, then we write $\text{Ker}\left(-\Delta_{g_{j}}-\lambda_{k}^{\left(j\right)}\right) = S_{k}^{(j)}\oplus \left(S_{k}^{(j)}\right)^{\perp}$,
       	   		and find a nontrivial $\phi_{k,l_{0}}^{(j)}\in \left(S_{k}^{(j)}\right)^{\perp}$.
       	   		Thus we have for any $u_{j}^{f}\in K^{(j)}$ in \eqref{DSS},
       	   		\begin{equation*}
       	   			0 = \left(\pi_{k}^{\left(j\right)}u_{j}^{f}, \phi_{k,l_{0}}^{\left(j\right)}\right)_{L^{2}\left(M\right)} = \left(u_{j}^{f}, \phi_{k,l_{0}}^{\left(j\right)}\right)_{L^{2}\left(M\right)},
       	   		\end{equation*}
       	   		which contradicts Lemma \ref{SK}.
       	   		   	   		
       	   		Combining \eqref{DR3} with \eqref{SDI}, we deduce that there exists a $L^{2}\left(M\right)$-basis  $\left\{\{\phi_{k,l}^{\left(j\right)}\}_{l=1}^{d_{k}^{(j)}}\right\}_{k=0}^{\infty}$ consisting of eigenfunctions of $-\Delta_{g_{j}}$ corresponding to eigenvalues 
       	   		$$0=\lambda_{0}^{\left(j\right)}<\lambda_{1}^{\left(j\right)}<\cdots<\lambda_{k}^{\left(j\right)}<\cdots$$
       	   		such that 
       	   		\begin{equation}\label{SRS}
       	   			\lambda_{k}^{\left(1\right)}=\lambda_{k}^{\left(2\right)}:=\lambda_{k},\: d_{k}^{\left(1\right)}=d_{k}^{\left(2\right)}:=d_{k}\quad \text{and}\quad \phi_{k,l}^{\left(1\right)}|_{O}=\phi_{k,l}^{\left(2\right)}|_{O},\quad k\in \mathbb{N},\, l=1,\cdots, d_{k}.
       	   		\end{equation}
       	   	This completes the proof of the equality in \eqref{SDI2}.
       	   	
Moreover, the orthogonality of $\left\{\{\phi_{k,l}^{\left(j\right)}\}_{l=1}^{d_{k}^{\left(j\right)}}\right\}_{k=0}^{\infty}$ for $j=1,2$ was proved under the condition \eqref{eqn : FKU condition} in \cite[Theorem 1.11]{feizmohammadi2024calderonproblemfractionalschrodinger}. 
       	\end{proof}
  
  \begin{proof}[Proof of Theorem \ref{MR}]
  	In summary, we determine the heat semigroup data \eqref{RM1I} with $H_{g,0,V}^{O}$ replaced by $H_{g,b,V}^{O}$ from the local source-to-solution map $\Lambda_{ g,b,V}^{O}$ in Proposition \ref{P1} and then recover the spectral data \eqref{SPD} in Proposition \ref{P2} together with Proposition \ref{P3}. Thus there is a diffeomorphism $\Psi:M\rightarrow M$ such that $$
  		g_{1}=\Psi^{*}g_{2}\quad \text{and} \quad \Psi|_{O}=\mathrm{Id}$$ from Lemma \ref{ML}.
  		With the gauge invariance of the local source-to-solution map in Theorem \ref{DG} and the uniqueness of coefficients on a priori known Riemannian manifold in Theorem \ref{MR1}, we obtain that $$b_{1}=\Psi^{*}b_{2}\quad \text{and} \quad V_{1}=\Psi^{*}V_{2}.$$
  		Thus we reconstruct all coefficients simultaneously and uniquely under a gauge class. 
  	\end{proof}

  	\bigskip
  
  \noindent {\bf Acknowledgements.} The authors were supported in part by NSFC grant 12571451. The authors thank Ali Feizmohammadi, Katya Krupchyk and Gunther Uhlmann for helpful comments and suggestions.

   \bigskip	\noindent {\bf Data Availability Statement.} Data sharing not applicable to this article as no datasets were generated or analysed during the current study.
  
   \bigskip	\noindent {\bf Conflict of Interest.} The authors have no conflicts of interest to declare that are relevant to the content of this article.

\bibliographystyle{abbrv}
\bibliography{ref}
\end{document}